\documentclass[natbib]{svjour3}

\RequirePackage{fix-cm}
\usepackage{mathptmx}  
\smartqed  
\usepackage{amsfonts,amsmath,amssymb,mathrsfs,epsfig,a4}
\usepackage{latexsym}
\usepackage{xcolor}
\usepackage{longtable}
\usepackage{multirow}

\usepackage{amsthm}

\allowdisplaybreaks
\spnewtheorem*{RGproof}{Proof:}{\itshape}{\rmfamily}

\newcommand{\bM}{\mathbb{M}}
\newcommand{\bN}{\mathbb{N}}

\newcommand{\bR}{\mathbb{R}}
\newcommand{\bC}{\mathbb{C}}
\newcommand{\bS}{\mathbb{S}}

\newcommand{\bB}{\mathbb{B}}
\newcommand{\bY}{\mathbb{Y}}

\newcommand{\cB}{\mathcal{B}}

\newcommand{\cF}{\mathcal{F}}
\newcommand{\cG}{\mathcal{G}}
\newcommand{\cH}{\mathcal{H}}

\newcommand{\cL}{\mathcal{L}}

\newcommand{\unif}{\text{\rm unif}}

\newcommand{\cov}{\text{\rm cov}}

\newcommand{\todistr}{\to_{\text{\rm\tiny distr}}} 
\newcommand{\eqdistr}{=_{\text{\rm\tiny distr}}}

\newcommand{\rk}{{\text{\rm rk}}}
\newcommand{\AV}{{\text{\rm AV}}}

\newcommand{\GL}{{\text{\rm GL}}}
\DeclareMathOperator{\Tr}{Tr} 

\allowdisplaybreaks

\newcommand{\PiF}{\Pi^{\tiny \rm F}}
\newcommand{\PiH}{\Pi^{\tiny \rm H}}
\newcommand{\LaP}{\Lambda^{\tiny \rm P}}
\newcommand{\LaC}{\Lambda^{\tiny \rm c}}
\newcommand{\PiC}{\Pi^{\tiny \rm c}}
\newcommand{\YT}{\bY_{\text{\tiny T}}}

\newcommand{\poF}{<_{\text{\tiny{F}}}}
\newcommand{\poH}{<_{\text{\tiny{H}}}}
\newcommand{\pH}{p^{\tiny \rm H}}

\journalname{}

\makeatletter
\def\LT@makecaption#1#2#3{%
  \LT@mcol\LT@cols c{\hbox to\z@{\hss\parbox[t]\LTcapwidth{%
    \captionstyle
    \sbox\@tempboxa{{\floatlegendstyle#1{#2}\floatcounterend}\capstrut #3}%
    \ifdim\wd\@tempboxa>\hsize
      {\floatlegendstyle#1{#2}\floatcounterend}\capstrut #3\par%
    \else
      \hbox to\hsize{\leftlegendglue\box\@tempboxa\hfil}%
    \fi
    \endgraf\vskip\baselineskip}%
  \hss}}}
\LTcapwidth=\dimexpr\textwidth-2\fboxsep\relax
\makeatother  

\begin{document}

\title{Ranks, copulas, and permutons}

\author{R. Gr\"ubel}

\institute{R. Gr\"ubel \at
           Institute of Actuarial and Financial Mathematics \\ 
 	       Leibniz Universit\"at Hannover\\
 	       Postfach 6009\\
           D-30060 Hannover, Germany\\
 	   \email{rgrubel@stochastik.uni-hannover.de}}

\date{Received: date / Accepted: date}

\maketitle

\begin{abstract}
We review a recent development at the interface between discrete mathematics on one hand
and probability theory and statistics on the other, specifically the use of Markov chains and
their boundary theory in connection with the asymptotics of randomly growing permutations. 
Permutations connect total orders on a finite set, which leads to the use of a pattern
frequencies. This view is  closely related to classical concepts of nonparametric statistics. 
We give several applications and discuss related topics and research areas, 
in particular the treatment of other combinatorial families, the cycle view of
permutations, and an approach via exchangeability. 

\end{abstract}

\section{Introduction}\label{sec:intro}

Ranks are at the basis of nonparametric statistics, copulas are a standard tool for
the description of dependencies, and permutons have recently been introduced as limit
objects for large permutations. We survey some of the connections between these concepts,
mainly from a probabilistic and statistical point of view. In particular, 
we explain the use of Markov chains and their boundary theory in the context 
of sequences of permutations. 
This approach has received some attention in connection with the asymptotics 
of other randomly growing discrete structures, such as graphs and trees.

The paper addresses a statistically educated audience. Consequently, we only 
mention Chapter 13 of~\cite{vdV} in connection with ranks in statistics, and for copulas
we refer to~\cite{Nelson}. In Section~\ref{sec:perm}
we recall some basic aspects of permutations; see \cite{Bona} for a textbook reference.
In Section~\ref{sec:subs} we introduce permutons, which are two-dimensional copulas.
The basic reference here is the seminal paper~\cite{Hopp}. The construction of limits is 
discussed in some generality, again from a probabilistic point of view. 

Permutons arise 
as limits if pattern frequencies are considered.
In Section~\ref{sec:dyn} we explain an approach to the
augmentation of discrete spaces that is based on the construction of suitable Markov 
chains and the associated boundary theory. The latter can be seen as a discrete variant 
of classical potential theory, going back to~\cite{Doob}; see \cite{WoessDMC} for a 
general textbook reference and~\cite{GrSemBerKMK} for an elementary
introduction closer to the present context. 
We show that the pattern frequency topology arises 
naturally if a specific Markov chain $(\Pi_n)_{n\in\bN}$ is chosen. 
This approach also provides an almost
sure construction, meaning that $\Pi_n\to\Pi_\infty$ `pathwise' with probability one. 

Section~\ref{sec:app} 
collects some applications. Finally, in the last section,
we sketch some related problems and research areas, where we use the reviewish character of
the present paper as a license for explaining connections that seem to be known to researchers in 
the field, but not necessarily to an interested newcomer (which the author was, some years 
ago). It should be clear that we will have to (and do) resist the temptation of always giving
the most general results.

Throughout, a specific data set will  be used to illustrate the various concepts,
and we occasionally point out analogies to statistical concepts.

\section{Permutations}\label{sec:perm}
For a finite set $A$ we write $\bS(A)$ for the set of bijections $\pi:A\to A$ and abbreviate 
this to $\bS_n$ if $A=[n]:=\{1,\ldots,n\}$. An element $\pi$ of $\bS_n$ can be represented 
in one-line notation or in standard cycle notation.
For the first $\pi=(\pi(1),\pi(2),\ldots,\pi(n))$ (or $\pi(1)\pi(2)\cdots\pi(n)$ if no
ambiguities can arise) is simply the tuple of its values, and 
we note at this early stage that this notation presupposes a total (or linear)  
order on the base set $A$.
The second one is somewhat closer to the notion of a bijection: Each element $a$ of $A$ 
will run through a finite cycle if $\pi\in\bS(A)$ is applied repeatedly. In the special case
$A=[n]$ this leads to a tuple of cycles, and standard notation means that these are 
individually ordered by putting the respective largest element first and then listing the cycles 
increasingly according to their first elements. The \emph{transition lemma}, also known as
\emph{Foata's correspondence}, provides a bijection between the two representations by  
erasing the brackets in the step from cycle to one-line notation, and, roughly, 
by inserting brackets in front of every increasing record in the one-line notation, 
such as in 
\begin{equation}\label{eq:Foata}
 \underline{5}\,\underline{8}7324\underline{9}16 \quad 
          \longrightarrow \quad       (5)(87324)(916). 
\end{equation}
Much of the classical material connecting permutations and probabilities refers to
the cycle structure, typically under the assumption that the random permutation $\Pi_n$
is chosen uniformly from $\bS_n$, which we abbreviate to $\Pi_n\sim\unif(\bS_n)$.
A standard problem, often discussed in elementary probability courses, 
concerns the probability $p_n$
that $\Pi_n$ has at least one fixed point, i.e.~a cycle of length $1$, which can be evaluated
with the inclusion-exclusion formula to $p_n=\sum_{k=1}^n (-1)^{k+1}/k!$. 
Further, the correspondence
in~\eqref{eq:Foata} relates the number of cycles on the right to the number of increasing
records on the left hand side. As a final example we mention that 
the distribution of the cycle
type $C_n=(C_n(1),C_n(2),\ldots)$ of $\Pi_n$, where $C_n(i)$ denotes the number of cycles
of length $i$ in $\Pi_n$, is given by
\begin{equation}\label{eq:Cauchy}
  P\bigl(C_n(1)=c_1,C_n(2)=c_2,\ldots\bigr)\; =\; \prod_{i=1}^\infty \frac{1}{i^{c_i}c_i!},
\end{equation}   
where the values $c_i$, $i\in\bN$, have to satisfy the obvious constraint
that $\sum_{i\in\bN} ic_i=n$. In particular, $C_n(i)=0$ for $i>n$. 
Multiplication by $n!$ leads to a formula for the number of permutations of a 
given cycle type. From a probabilistic point of view, rewriting~\eqref{eq:Cauchy} as
\begin{equation}\label{eq:cycrepr}
  \cL\bigl((C_n(1),\ldots,C_n(n))\bigr) \; =\; \cL\bigl[(Z_1,\ldots,Z_n)\big|T_n=n\bigr],
\end{equation}
with $Z_1,\ldots,Z_n$ independent, $Z_i$ Poisson distributed
with parameter $1/i$, and $T_n:=\sum_{i=1}^n i\, Z_i$, seems more instructive; for example, 
it provides access to the distributional asymptotics. We refer
to~\cite{ABT} for this and related
material. Our later aim, however, is to deal with limits for the permutations 
themselves, but we will return to the cycle view in Sections~\ref{subsec:cycle}.

We may also describe a permutation $\pi\in\bS_n$ by its $n\times n$ permutation matrix $M(\pi)$ with entries $m_{ij}=1$ if $\pi(j)=i$ and $m_{ij}=0$ if $\pi(j)\not= i$. For the permutation on the left in~\eqref{eq:Foata} we get
\begin{equation}\label{eq:permmatrix}  
M(\pi)\, =\, 
\begin{pmatrix} 
           0&0&0&0&0&0&0&1&0 \\ 0&0&0&0&1&0&0&0&0\\ 0&0&0&1&0&0&0&0&0\\
           0&0&0&0&0&1&0&0&0 \\ 1&0&0&0&0&0&0&0&0\\ 0&0&0&0&0&0&0&0&1\\
           0&0&1&0&0&0&0&0&0 \\ 0&1&0&0&0&0&0&0&0\\ 0&0&0&0&0&0&1&0&0
\end{pmatrix} .
\end{equation}
Below we will return to the fact that the composition of permutations corresponds to 
matrix multiplication in this description.

Finally,  permutations appear as order isomorphisms: If $<_1$ and $<_2$ are total orders 
on some finite set $A$, then there is a unique $\pi\in\bS(A)$ such that
\begin{equation*}
  a<_1 b\  \Longleftrightarrow\ \pi(a) <_2 \pi(b)\qquad\text{for all } a,b\in A.
\end{equation*}
A one-line notation such as $587324916$ in \eqref{eq:Foata}
defines a total order on the set of numerals, with $5$ the smallest and $6$ 
the largest element. With `$<$' the canonical order on $[n]$ (or $\bR$) 
the permutation $\pi\in\bS_n$ 
then provides an order isomorphism between $([n],<)$ and $([n],<_\pi)$, with `$<_\pi$' 
given by the one-line notation for $\pi$. 

Suppose now that we have data $x_1,\ldots,x_n\in\bR$, with $x_i\not=x_j$ if $i\not=j$  (no ties).
The \emph{rank} of $x_i$ in this set is given by 
$\rk(x_i)=\rk(x_i|x_1,\ldots,x_n)=\#\{j\in [n]:\, x_j\le x_i\}$,
and these can be combined into a permutation 
$\pi_x:=(\rk(1),\ldots,\rk(n))\in\bS_n$. 
For two-dimensional data $(x_1,y_1),\ldots,(x_n,y_n)$ with no ties in the two sets of 
component values 
we then obtain two permutations $\pi_x,\pi_y$ that define two total orders on $[n]$. 
These are, in the above sense, connected by the permutation
\begin{equation}\label{eq:connect}
        \pi=\pi((x_1,y_1),\ldots,(x_n,y_n))=\pi_y\circ\pi_x^{-1}.
\end{equation}
In this situation, the permutation plot for $\pi$ is known as the  
\emph{rank plot} for the data. The order-relating 
permutation does not depend on the order in which the observation vectors are numbered.
Formally, it is invariant  under permutations $\sigma$ of $[n]$
in the sense that the permuted pairs 
$(x_{\sigma(1)},y_{\sigma(1)}),\ldots,(x_{\sigma(n)},y_{\sigma(n)})$
lead to the same $\pi$.
Hence, in statistical terms, $\pi$ depends on the data only through 
their empirical distribution.

\begin{figure}
\setlength{\abovecaptionskip}{-.25cm}

\hbox{{}\hspace{-1.45cm}\includegraphics[scale=.5]{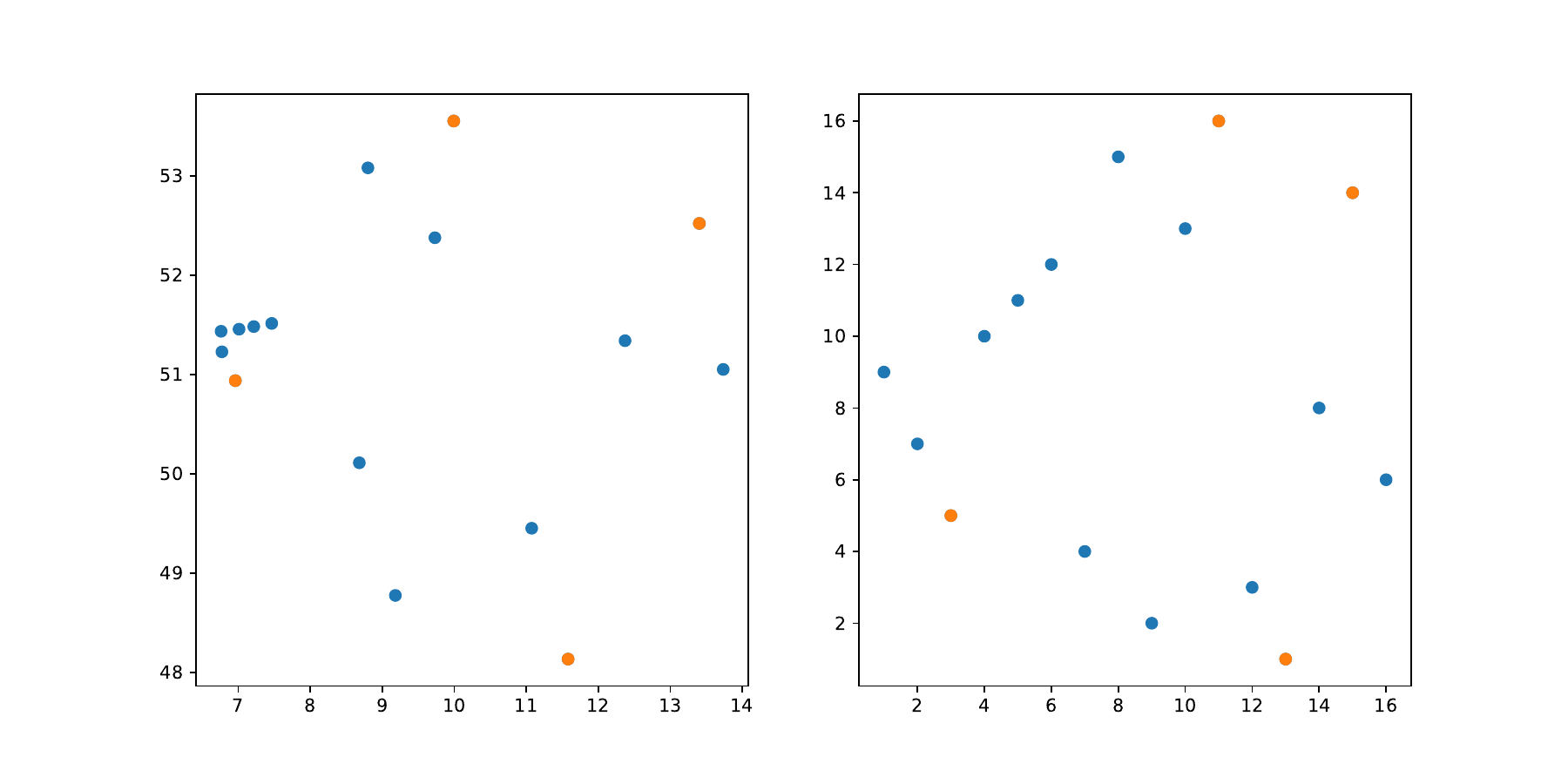}}

\caption{Scatter plot and rank plot for the city data (red: the four largest cities)}
\label{fig:Staedte}
\end{figure}

We close this section with a data example. The 16 largest cities in Germany may be ordered 
alphabetically, or by decreasing population size, or by location in the west-east or south-north 
direction. Starting with an ordering by size the corresponding geographic orderings 
are
\begin{gather*}
   \pi_{\tiny \text{long}} =(15, 11, 13,  3,  7,  9,  2,  6,  4,  8, 16, 14, 10, 12,  1,5),\\
   \pi_{\tiny \text{lat}} = (14, 16,  1,  5,  4,  2,  7, 12, 10, 15,  6,  8, 13,  3,  9, 11).
\end{gather*}
For example, Berlin is fairly far in the east and fairly far up north within
this group. The permutation relating these as total orders is given by
\begin{equation}\label{eq:rank}
     \pi_{\tiny \text{lat}} \circ \pi_{\tiny \text{long}}^{-1} 
           = (9,  7, 5,  10,  11, 12, 4,  15,  2, 13, 16,  3, 1, 8, 14,  6).
\end{equation}
Figure~\ref{fig:Staedte} shows the true locations (in longitude and latitude degrees)
on the left and the rank plot on the right hand side, with the four largest cities in red. 
The observant reader will have noticed that the rank plot is, essentially, the associated
permutation matrix. In fact, the plot of a permutation also leads to an interpretation
(or coding) of permutations 
that will be central for the permuton aspect: Writing $\delta_x$ for the
one-point measure at $x$ we can associate a distribution $\mu(\pi)$ on the unit square
with $\pi\in\bS_n$ via
\begin{equation}\label{eq:discrCop}
     \mu(\pi) \;:=\, \frac{1}{n} \sum_{i=1}^n \delta_{(i/n,\pi(i)/n)}.
\end{equation} 
Clearly, $\pi$ can be recovered from $\mu(\pi)$.
Also, for all $\pi\in\bS_n$, the marginals of $\mu(\pi)$ are the discrete 
uniform distributions on the set $\{1/n,2/n,\ldots,1\}$.

In view of later developments we note that in the passage from rank plot to permutation
plot the relative positions of the points in $x$- and $y$-direction do not change. 
In our data example,  a city that is north (up) or east (right) of some other city will 
stay so if we move from the left to the right part of~Figure~\ref{fig:Staedte},

\section{Patterns, subsampling and convergence}\label{sec:subs}

Our aim in this section is to obtain a formal framework for the informal 
(and ubiquitous) question, 
`what happens as $n\to\infty$?'. We may, for example, 
start with a sequence $(\pi_n)_{n\in\bN}$ of permutations with $\pi_n\in\bS_n$ for all
$n\in\bN$ and regard these as elements of the set   
$\bS:=\bigsqcup_{n=1}^\infty \bS_n$, which we then need to augment in order to obtain 
a limit. A standard procedure begins with a set $\cF$ of functions $f:\bS\to [0,1]$ 
that separates the points in the sense that, for each pair $\pi,\sigma\in\bS$ with
$\pi\not=\sigma$, there is a function $f\in\cF$ such that $f(\pi)\not= f(\sigma)$. 
We then identify $\pi\in\bS$ with the function $f\mapsto f(\pi)$, which is an element 
of the space $[0,1]^{\cF}$ which, endowed with the product topology, is compact by 
Tychonov's theorem. 
The closure of the embedded $\bS$ provides a compactification $\bar\bS$ of $\bS$,
where convergence of a sequence $(\pi_n)_{n\in\bN}$ is equivalent to the convergence 
of all real sequences $(f(\pi_n))_{n\in\bN}$, $f\in\cF$. Moreover, 
each function $f\in\cF$ has a continuous extension to $\bar \bS$.
In the present situation $\cF$ will always be countable, and we start with the 
discrete topology on $\bS$. Then $\bar\bS$ is a separable compact topological space,
and the boundary $\partial \bS$ is a compact subset. 

We discuss three function classes and the resulting topologies. 

First, let $\cF$ be the set of indicator functions $f_\sigma=1_{\{\sigma\}}$, $\sigma\in\bS$,
i.e.~$f_\sigma(\sigma)=1$ and $f_\sigma(\pi)=0$ if $\pi\not=\sigma$. Then, if  
$(\pi_n)_{n\in\bN}$ is a sequence that leaves $\bS$ in the sense that for each $k\in\bN$ there 
is an $n_0=n_0(k)\in\bN$ such that $\pi_n\notin\bS_k$ for all $n\ge n_0$, we always get the
function which is equal to 0 as the pointwise limit, hence the boundary consists of a single 
point. This topology is known as the \emph{one-point compactification}.

For the second example we note first that order relations transfer to subsets by restriction. 
For $A=[n]$, $\pi\in\bS_n$, and $B=[k]$ with $k\le n$, restricting $<_\pi$ from $A$ to $B$
amounts to deleting all entries greater than $k$ in the one-line notation of $\pi$, which
gives the one-line notation for a permutation $\sigma\in\bS_k$. We generally write 
$\sigma=\pi[B]$ if $(B,<_\pi)$ is order isomorphic to $([k],<_\sigma)$. If 
$B=\{i_1,\ldots,i_k\}\subset [n]$ with $i_1<\cdots <i_k$ this is equivalent to 
\begin{equation*}
 \sigma(j)<\sigma(m) \ \Longleftrightarrow\ \pi(i_j)<\pi(i_m)\quad \text{for } j,m\in [k],
\end{equation*}
a condition that will appear repeatedly below.
We now define $f_\sigma:\bS \to [0,1]$, $\sigma\in\bS_k$, 
by $f_\sigma(\pi)=1$ if $\sigma=\pi[A]$ for $A=[k]$, 
and $f_\sigma(\pi)=0$ otherwise. As in the previous case the functions take only the 
values 0 and 1, so that convergence means that $f_\sigma(\pi_n)$ remains constant from 
some $n_0=n_0(\sigma)$ onwards. 
In view of $\sum_{\sigma\in\bS_k} f_\sigma(\pi)= 1$ for all $\pi\in\bS_n$,
$n\ge k$, any such sequence leads to a sequence $(\sigma_k)_{k\in\bN}$, $\sigma_k\in\bS_k$,
that has the property that $f_{\sigma_k}(\pi_n)=1$ for all $n\ge n_0(k)$. Clearly, this sequence
must be projective in the sense that the restriction of $\sigma_k$ to $[l]$ with $l\le k$ is
equal to $\sigma_l$. This implies that there is a total order $<_{\pi_\infty}$ on $\bN$ which 
reduces to $<_{\sigma_k}$ on $[k]$ for all $k\in\bN$. The result is known as the 
\emph{projective topology}, and the set of limits is the set of total orders on~$\bN$.

The third topology will lead to permutons. It arises if we use \emph{pattern counting}, 
another aspect 
that is closely related to the interpretation of permutations as connecting total orders.
Let $M(k,n)$ be the set of strictly increasing functions $f:[k]\to[n]$, $n\ge k$. We
define the \emph{relative frequency of the pattern $\sigma$ in the permutation $\pi$} as
 \begin{equation*}
 t(\sigma,\pi) \, =\, \frac{1}{\# M(k,n)}\,\#\bigl\{f\in M(k,n):\, 
                    \sigma(i)<\sigma(j)\, \Leftrightarrow\, \pi(f(i))<\pi(f(j)) \bigr\}.
\end{equation*}
Of course, $\# M(k,n)=\binom{n}{k}$. We augment this by putting $t(\sigma,\pi)=0$ if $k> n$.
Using the notation from the previous paragraph and identifying $f\in M(n,k)$ with its 
range $A$ we see that the numerator is the number  
of size $k$ subsets $A$ of $[n]$ with $\pi[A]=\sigma$.   
Table~\ref{tab:staedte} contains the absolute and relative frequencies for the six patterns 
of length three in the city data introduced at the end of Section~\ref{sec:intro}. The four
largest cities generate the pattern $\sigma=2413\in\bS_4$.

\begin{table}\label{tab:staedte}
\begin{tabular}{|l|rrrrrr|}
\hline 
\rule[-1ex]{0pt}{3.5ex}permutation &123 & 132 & 213 & 231 & 312 & 321 \\ 
pattern frequency &69 & 57 & 88 & 61 & 130 & 50 \\ 
relative frequency &.152 & .125 & .193 & .134 & .286 & .110 \\ 
\hline 
\end{tabular} 
\vspace{2mm}
\caption{Patterns in the city data permutation}
\end{table}


From a probabilistic or, in fact, statistical point of view, we may regard the
relative frequencies as probabilities in a sampling experiment. With $\xi=\xi_{n,k}$ 
uniformly distributed on the set of all subsets $A$ of $[n]$ with size $k$, 
the function $\sigma\mapsto t(\sigma,\pi)$ 
is the probability mass function of the random element $\pi[\xi]$ of $\bS_k$.
In particular, whenever $k\le n$, 
\begin{equation}\label{eq:pmf}
     \sum_{\sigma\in\bS_k}  t(\sigma,\pi) \, = \, 1\quad \text{for all } \pi\in\bS_n.
\end{equation}
Building on an earlier approach in graph theory, \cite{Hopp} introduced a notion 
of convergence based on such pattern frequencies. 
We identify a permutation $\pi$ with its function $\sigma\to t(\sigma,\pi)$ on the
set $\bS$ and with values in the unit interval; it is easy to see that different
permutations lead to different functions.  
The general approach outlined above then leads us to define 
convergence of the sequence $(\pi_n)_{n\in\bN}$ as convergence of 
$t(\sigma,\pi_n)$ for all $\sigma\in\bS$.
Below we will refer to this as the \emph{pattern frequency topology},

For the construction to be of use we need a 
description of the  boundary $\partial\bS=\bar\bS\setminus \bS$.
In the general setup the limit of a convergent sequence $(\pi_n)_{n\in\bN}$ with 
$|\pi_n|\to\infty$ is simply the function 
\begin{equation*}
 \pi_\infty:\bS\to [0,1], \quad \sigma\mapsto \lim_{n\to\infty}t(\sigma,\pi_n),
\end{equation*}
but obviously, not all functions from $\bS$ to $[0,1]$ 
appear in that manner. For example, \eqref{eq:pmf} implies that
$\sum_{\sigma\in\bS_k}\pi_\infty(\sigma)=1$ for any fixed $k\in\bN$, and a conditioning argument 
leads to the inequality $\pi_\infty(\sigma)\ge  t(\sigma,\rho)\cdot \pi_\infty(\rho)$ for
all $\sigma\in\bS_k$, $\rho\in\bS_n$ with $k\le n$. In fact, we even have, for all $k,l,m\in\bN$
with $k<l<m$ and all $\sigma\in\bS_k$, $\rho\in\bS_m$,
\begin{equation}\label{eq:cotranspattern}
    t(\sigma,\rho)\, =\, \sum_{\tau\in\bS_l} t(\sigma,\tau)\, t(\tau,\rho).
\end{equation} 
This follows from viewing the sampling procedure as removing single elements repeatedly,
and then decomposing with respect to the permutation obtained at the intermediate level. 
The decomposition \eqref{eq:cotranspattern} will play an important role later in 
the context of the dynamics of permutation sequences; see Section~\ref{sec:dyn}.

Now copulas enter the stage: A (two-dimensional) \emph{copula} $C:[0,1]\times[0,1]\to[0,1]$ 
is the distribution function of a probability measure on the unit square that has uniform
marginals, meaning  that $C(u,1)=u$ for all $u\in [0,1]$ and $C(1,v)=v$ for all
$v\in[0,1]$. We define the associated \emph{pattern frequency function} $t(\cdot,C)$ by sampling:
For $\sigma\in\bS_k$ let $t(\sigma,C)$ be the probability that a sample 
$(X_1,Y_1),\ldots,(X_k,Y_k)$ from $C$ leads to the permutation $\sigma$ in the sense that
\begin{equation*}
   \Pi_n((X_1,Y_1),\ldots,(X_k,Y_k)) = \sigma, 
\end{equation*}
see also~\eqref{eq:connect}. 
As the joint distribution of the sample is invariant under permutations it follows that,
for all $k\in\bN$ and $\sigma\in\bS_k$,
\begin{equation*}
   t(\sigma,C) \;=\; k!\, P(X_1<X_2<\cdots < X_k,\, 
              Y_{\sigma(1)}<Y_{\sigma(2)}<\cdots< Y_{\sigma(k)}). 
\end{equation*}
We observe (but do not prove) that $C$ is determined by the function $\sigma\mapsto t(\sigma,C)$
on $\bS$, and that for copulas, weak convergence is equivalent to the pointwise
convergence on $\bS$ of these functions. Note the analogy to the moment problem
in the first and to the moment method in the second statement.   

We can now collect and rephrase some of the main results in~\cite{Hopp}. For the 
definition of the random permutations in part (c) we refer to~\eqref{eq:connect} again, and to
the description given there  of a suitable algorithm for obtaining the permutation from the data.
Ties in the $X$- resp.~$Y$-values may be ignored as these are, individually, a sample
from the uniform distribution on the unit interval.

\begin{theorem}\label{thm:Hoppen}
\emph{(a) } A sequence $(\pi_n)_{n\in\bN}$ of permutations in $\bS$ 
with $|\pi_n|\to\infty$ 
converges in the pattern frequency topology if and only if there is a copula $C$
such that 
\begin{equation}\label{eq:limC}
      \lim_{n\to\infty}t(\sigma,\pi_n)=t(\sigma,C)\quad\text{for all } \sigma\in\bS.
\end{equation}
Moreover, the copula $C$ is unique.

\vspace{.5mm}
\emph{(b) } In the pattern frequency topology, the boundary $\partial\bS$ is
homeomorphic to the space of copulas, endowed with weak convergence.

\vspace{.5mm}
\emph{(c) } Let $C$ be a copula and suppose that $(X_i,Y_i)$, $i\in\bN$, are independent 
random vectors with distribution function $C$. For each $n\in\bN$ let $\Pi_n$ be the random
permutation associated with the first $n$ pairs $(X_1,Y_1),\ldots,(X_n,Y_n)$.
Then, in the pattern frequency topology, $\Pi_n$ converges to $C$ almost surely
as $n\to\infty$.
\end{theorem}

Below we will refer the stochastic process $(\Pi_n)_{n\in\bN}$ in part (c) 
of Theorem~\ref{thm:Hoppen} as 
\emph{the permutation sequence generated (or parametrized) by the copula $C$}.

From a probabilistic point of view it is not too difficult to understand this result:
Let $\bM_1=\bM_1([0,1]\times [0,1], \cB([0,1]\times [0,1]))$
be the space of probability measures on the Borel subsets of the unit square, endowed
with the topology of weak convergence, and consider the embedding $\pi\mapsto \mu(\pi)$
of $\bS$ into $\bM_1$ given by~\eqref{eq:discrCop}. Suppose now that 
$(\pi_n)_{n\in\bN}$ is a sequence of permutations where we assume for simplicity that 
$\pi_n\in\bS_n$ for all $n\in\bN$. As $\bM_1$ is compact with respect to weak convergence
this sequence has a limit point $\mu_\infty$, and convergence would follow from the
uniqueness of $\mu_\infty$. For this we first note that the marginals of any such limit
point are the uniform distributions on the unit interval, as weak convergence implies the
convergence of the marginal distributions. Hence any limit point 
$\mu_\infty$ has as its distribution function a copula $C$. 
Next we argue that~\eqref{eq:limC} holds: Let $\sigma\in\bS_k$ be given and let the 
function $h_\sigma$
on $([0,1]\times [0,1])^k$ indicate whether or not the order-relating permutation
for the marginal rank vectors of $((x_1,y_1),\ldots,(x_k,y_k))$ is equal to $\sigma$. Then
\begin{equation}\label{eq:Urepr}
       t(\sigma,\pi_n) 
           \,=\,     \frac{1}{\binom{n}{k}} \sum_{1\le i_1<\cdots<i_k\le n}
         h_\sigma\bigl((x_{i_1},y_{i_1})\ldots,(x_{i_k},y_{i_k})\bigr). 
\end{equation}
Uniform marginals implies that the set of discontinuities of $h_\sigma$ is a  
$\mu_\infty$-null set, and the weak convergence of $\mu(\pi_n)$ to $\mu_\infty$
implies the weak convergence of the respective $k$th measure-theoretic powers.
Taken together this shows~\eqref{eq:limC}. The desired uniqueness now follows
from the fact that a copula $C$ is determined by its function $\sigma\to t(\sigma,C)$.
For (b) we use that weak convergence of copulas 
is equivalent to pointwise convergence of the associated functions on $\bS$. 
Finally, part (c) is immediate from the representation~\eqref{eq:Urepr} of pattern 
frequencies as $U$-statistics, and the classical
consistency result for the latter.

\section{Dynamics}\label{sec:dyn}

Generating random permutations is an important applied task. 
In line with the 
history of probability we recall a standard gaming example: Card mixing requires 
$\Pi_n\sim  \unif(\bS_n)$ for some fixed $n$. Many shuffling algorithms may be 
modeled as
random walks on $\bS_n$ and can then be analyzed using the representation
theory for this group; see also Section~\ref{subsec:cycle} below. 

Here we are interested in sequences $(\Pi_n)_{n\in\bN}$ with $\Pi_n\in\bS_n$ for all
$n\in\bN$, and a Markovian dynamic. A classical example is the 
Markov chain  $\PiF=(\PiF_n)_{n\in\bN}$ which starts
at $\PiF_1\equiv (1)$ and then uses independent random variables $I_n\sim\unif([n+1])$, 
independent of $\PiF_1,\ldots,\PiF_n$, 
to construct $\PiF_{n+1}$ from $\PiF_n=(\pi(1),\ldots,\pi(n))$ as follows,
\begin{equation}\label{eq:transF}
  \Pi^F_{n+1}=\bigl(\pi(1),\ldots,\pi(I_n-1),n+1,\pi(I_n),\ldots,\pi(n)\bigr).
\end{equation} 
In words: We insert the new value $n+1$ at a position in a gap chosen uniformly at random.
This appears, for example, on p.132 in~\cite{FellerI}. 
A second possibility, which we will relate to~\cite{Hopp}, 
and which we denote by $\PiH=(\PiH_n)_{n\in\bN}$, uses `double insertion' instead.
This means that the step from $n$ to $n+1$ is based on two independent
random variables $I=I_n$ and $J=J_n$, both uniformly distributed on $[n+1]$, and
\begin{equation}\label{eq:transH}
      \PiH_{n+1}(i) \, =\, \begin{cases} \PiH_n(i), &\text{if } i<I\text{ and }\PiH_n(i)<J,\\
                          \PiH_n(i)+1, &\text{if }i<I\text{ and } \PiH_n(i)\ge J,\\
                           J, &\text{if } i=I,\\
                          \PiH_n(i-1), &\text{if } i>I\text{ and } \PiH_n(i)<J,\\
                          \PiH_n(i-1)+1, &\text{if }i>I\text{ and } \PiH_n(i)\ge J.
                           \end{cases}
\end{equation}
In words: We insert $J$ at position $I$ and increase the appropriate values by 1.
Obviously, this process is a Markov chain, again with start at the single element $\sigma=(1)$
of $\bS_1$. We collect some properties of this process. The notion of a permutation 
sequence generated by a copula is explained in Theorem~\ref{thm:Hoppen}\,(c).

\begin{lemma}\label{lem:Hoppen} 
Let $\PiH$ be as defined above. 

\vspace{0.6mm}
\emph{(a)} $\ \PiH$ is equal in distribution to the permutation sequence generated by 
the independence copula $C(x,y)= x\cdot y$.

\vspace{0.6mm}
\emph{(b)} The transition probabilities of $\PiH$ are given by
\begin{equation}\label{eq:transPiH}
      \pH(\sigma,\tau):= P(\PiH_{n}=\tau|\PiH_k=\sigma) \, =\, \frac{k!}{n!}\, t(\sigma,\tau)
\end{equation}
for all $k,n\in\bN$ with $k \le n$, and all $\sigma\in\bS_k$ and $\tau\in \bS_{n}$.

\vspace{1.2mm}
\emph{(c)} $\ \PiH_n\sim\unif(\bS_n)\ $ for all $n\in\bN$. 
\end{lemma}

\begin{RGproof}
(a) Let $\Pi=(\Pi_n)_{n\in\bN}$ be generated by the independence copula. 
Then the step from $\Pi_n$ to $\Pi_{n+1}$ is based on the component ranks $I$ and $J$ 
of the next pair $(X_{n+1},Y_{n+1})$
in the sequence $(X_1,Y_1),\ldots,(X_{n+1},Y_{n+1})$. It is well
known that $I$ and $J$ are both uniformly distributed on $[n+1]$. Here they are also
independent, hence $\Pi_{n+1}$ arises from $\Pi_n$ as in~\eqref{eq:transH}.

(b) The transition from $\sigma\in\bS_n$ to 
$\tau\in\bS_{n+1}$ is based on independent random variables $I$ and
$J$, both uniformly distributed on $[n+1]$. If $I=i$ then, with $\phi_i:[n]\to [n+1]$,
$\phi_i(k)=k$ if $k\le i$ and $\phi_i(k)=k+1$ if $k>i$, we must have
$\sigma(l)<\sigma(m)$ if and only if $\tau(\phi_i(l))<\tau(\phi_i(m))$ for all $l,m\in [n]$.
The number of such $i$'s is $(n+1)\,t(\sigma,\tau)$, and in each case we additionally need 
$J=\tau(I)$. This proves the one-step version of~\eqref{eq:transPiH}. 

To obtain the general case we now use induction, the Markov property,
and~\eqref{eq:cotranspattern}.

(c) This is immediate from \eqref{eq:transPiH} with $k=1$ and $\sigma=1$.
\qed
\end{RGproof}

\begin{figure}
\setlength{\abovecaptionskip}{-0.15cm}

\hbox{{}\hspace{-2.1cm}\includegraphics[scale=.46]{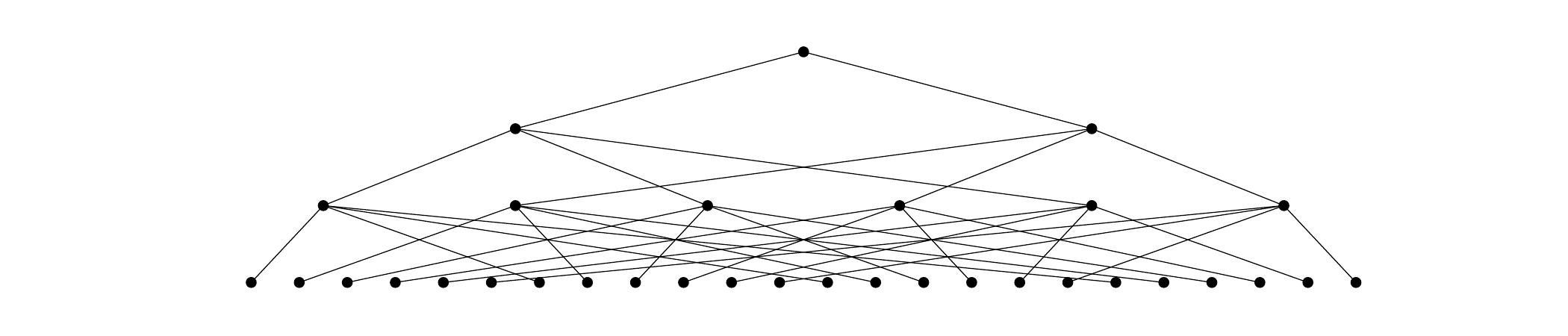}}

\vspace{-2mm}

\hbox{{}\hspace{-2.1cm}\includegraphics[scale=.46]{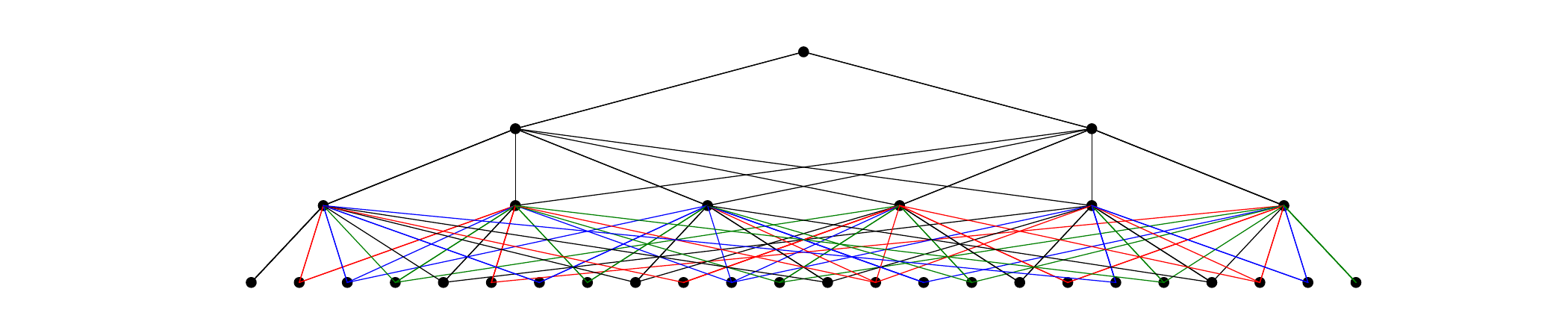}}

\caption{The transition graphs for the Markov chains $\PiF$ (top) and $\PiH$ (bottom)}
\label{fig:Hasse}
\end{figure}

Hence both $\PiF$ and $\PiH$ consist of uniformly distributed variables.
A trivial third possibility to achieve  this is to use independent $\Pi_n$'s with
$\Pi_n\sim\unif(\bS_n)$. The \emph{transition graph} for these chains has $\bS$ as its
set of nodes and an edge $\{\sigma,\tau\}$ for each pair $\sigma\in\bS_n$, $\tau\in\bS_{n+1}$,
$n\in\bN$, with the property that $P(\Pi_{n+1}=\tau|\Pi_n=\sigma)>0$.
Figure~\ref{fig:Hasse} shows the first four levels of the transition graphs for 
$\PiF$ and $\PiH$. These graphs are the \emph{Hasse diagrams} associated with two specific partial 
orders on $\bS$. In the first, $\sigma \poF \rho$ for $\sigma\in\bS_k$ 
and $\rho\in\bS_n$ with $k<n$ means that the one-line notation for $\sigma$ arises 
from the one-line notation for $\rho$ by simply deleting all values greater than $k$ 
in the latter, 
whereas in connection with $\PiH$ we have $\sigma \poH \rho$ if and only if
$\sigma$ appears as a pattern in $\rho$. 
Clearly, the graph for $\PiF$ is a tree: All previous states are deterministic functions
of the current state. Informally speaking, all three chains eventually leave 
the state space $\bS$, and this is what provides the connection to the previous section.

We need an excursion. The Markov property 
is often stated informally as `the future development depends on the history of the  process
only via its present state'. It may be instructive to carry out a simple calculation: If
$X_1,X_2,\ldots,X_n$ are random variables on some probability space with values in some
countable set, then their joint distribution can always be evaluated step by step via
\begin{equation*}
 P(X_1=x_1,\ldots,X_n=x_n) \;=\; P(X_1=x_1)\cdot 
                   \prod_{i=2}^n P(X_i=x_i|X_{i-1}=x_{i-1},\ldots,X_1=x_1),
\end{equation*}  
whenever the probability on the left hand side is positive. 
If the tuple has the Markov property then elementary manipulations of 
conditional probabilities lead to
\begin{align*}
 P(X_1=x_1,\ldots,X_n=x_n) \ &=\  P(X_1=x_1)\cdot \,
                   \prod_{i=2}^n P(X_i=x_i|X_{i-1}=x_{i-1}) \\
                            &=\  P(X_1=x_1)\cdot 
                   \prod_{i=2}^n P(X_{i-1}=x_{i-1}|X_i=x_i) 
                                    \,\frac{P(X_i=x_i)}{P(X_{i-1}=x_{i-1})}\\
             &=\ P(X_n=x_n) \prod_{i=n-1}^1 P(X_i=x_i|X_{i+1}=x_{i+1}).
\end{align*}  
Hence the Markov property then also holds `in the other direction'. Moreover, the converse
statement follows on reading the above from bottom to top. For a sequence $(X_n)_{n\in\bN}$,
and with $\cG_n,\cH_n$ the $\sigma$-fields generated  by 
the variables $X_1,\ldots,X_n$ and $X_n,X_{n+1},\ldots$ respectively, the 
Markov property can formally be written as 
$\cL(X_{n+1}|\cG_n)=\cL(X_{n+1}|X_n)$, $n\in\bN$, and the above argument shows that 
this implies $\cL(X_{n-1}|\cH_n)=\cL(X_{n-1}|X_n)$, $n\in\bN$, a procedure that may be 
interpreted as \emph{time reversal}. In fact, a somewhat neutral statement would be that, 
for a Markov process, past and future are conditionally independent, given the present. 

The distribution of a Markov
chain is specified by the distribution of the first variable and the (forward)
transition probabilities;
for the reverse chain we speak of the \emph{cotransition probabilities}. The tree structure 
of its transition graph implies that for $\PiF$ the cotransitions are degenerate in the sense
only the values 0 and 1 occur. For $\PiH$ we use Lemma~\ref{lem:Hoppen} to obtain
\begin{equation}\label{eq:cotransPiH}
   P(\PiH_n=\sigma|\PiH_{n+1}=\tau)\ 
             =\ \frac{P(\PiH_n=\sigma)}{P(\PiH_{n+1}=\tau)}\, \pH(\sigma,\tau)\
             =\  t(\sigma,\tau)
\end{equation}
for all $n\in\bN$, $\sigma\in\bS_n$ and $\tau\in\bS_{n+1}$. 
The decomposition~\eqref{eq:cotranspattern} can now be interpreted as describing the 
backwards movement of $\PiH$: \emph{ Pattern frequencies arise as cotransitions for this chain}.

Next we sketch Markov chain boundary theory, which grew out of the historically very 
important probabilistic approach to classical potential theory. For the highly transient 
chains considered 
here, where no state can be visited twice and where the 
time parameter is a function of the state and thus time homogeneity is implicit,
this takes on a particularly simple form. Suppose that $X=(X_n)_{n\in\bN}$ is a Markov chain 
with graded state space $S=\bigsqcup S_n$, that $X$ is adapted to the grading in the sense
that $X_n\in S_n$ for all $n\in\bN$, and that $S$ is minimal in the sense that $P(X_n=x)>0$
for all $x\in S_n$, $n\in\bN$. The \emph{Martin kernel} 
$K: S\times S\to \bR_+$ for $X$ is then given by
\begin{equation}\label{eq:Martin}
 K(x,y) \,:=\, \frac{P(X_n=y|X_m=x)}{P(X_n=y)}\quad \text{for all } x\in S_m,\, y\in S_n
\end{equation} 
if $n\ge m$, and $K(x,y)=0$ otherwise.
For each $x\in S_m$ let $f_x: S \to [0,1]$ be defined by 
\begin{equation}\label{eq:funcDM}
    f_x(y)\,=\,P(X_m=x)\, K(x,y)
           \, =\, P(X_m=x | X_n=y)
\end{equation} 
and let $\cF:=\{f_x:\, x\in S\}$. The \emph{Doob-Martin compactification}
$\bar S$ and the \emph{Martin boundary} $\partial S$  of the state space are the results of 
the procedure outlined at the beginning of Section~\ref{sec:subs}, with this choice of $\cF$ 
as space of separating functions. 

We give two central results.

\begin{theorem}\label{thm:harm}
In the Doob-Martin compactification, $X_n$ converges almost surely to a random variable
$X_\infty$ with values in the Martin boundary.
\end{theorem}

\begin{RGproof}
Let $f=f_x\in\cF$ with some $x\in S_k$. 
In view of the above remarks on the Markov property 
and time reversal,
\begin{equation*}
   f(X_n) \,=\,  P(X_k=x|X_n) \,=\, E[1_{\{x\}}(X_k)|\cH_n]\quad\text{for all }n\ge k,
\end{equation*}
so that, as $\cH_{n+1}\subset\cH_n$,
\begin{equation*}
   E[f(X_n)|\cH_{n+1}] \, =\, E\bigl[E[1_{\{ x\}}(X_k)|\cH_n]\big|\cH_{n+1}\bigr]
       \, =\, E[1_{\{x\}}(X_k)|\cH_{n+1}]\, =\, f(X_{n+1}).
\end{equation*}
This shows that $\bigl(f(X_n),\cH_n\bigr)_{n\ge k}$ is a reverse martingale, 
hence $f(X_n)$ converges almost surely. As $x\in S$ was arbitrary, the first statement 
follows by construction of the compactification.

As $P(X_\infty\in\bigcup_{k=1}^n S_k)=0$ for all $n\in\bN$, we also get 
$P(X_\infty\in\partial S)=1$.
\qed
\end{RGproof}

The general compactification procedure implies that the Martin kernel can be extended
continuously from  $S\times S$ to $S\times \bar S$ via $K(x,\alpha)=\lim_{n\in\bN} K(x,y_n)$
 where $(y_n)_{n\in\bN}$ is such that $y_n\to \alpha\in\partial S$ in the Doob-Martin topology.
The following is a special case of \emph{Doob's $h$-transform}. 

\begin{theorem}\label{thm:htrans} Let $X,S,\bar S,\partial S$ be as above and let 
$p(x,y)=P(X_{n+1}=y|X_n=x)$, $x\in S_n$, $y\in S_{n+1}$, be the transition function 
of the Markov chain $X$.
For $\alpha\in\partial S$ let $p^\alpha$ be defined by
\begin{equation*}
p^\alpha(x,y) \, =\,  \frac{K(y,\alpha)}{K(x,\alpha)}\, p(x,y)
                                   \quad\text{for all } n\in\bN,\; x\in S_n,\, y\in S_{n+1}.
\end{equation*}
Then $p_\alpha$ is the transition function for a Markov chain 
$X^\alpha=(X^\alpha_n)_{n\in\bN}$, and $X^\alpha_n$ converges almost surely to $\alpha$
in the Doob-Martin topology.
\end{theorem}

The picture that emerges from these results is that of an \emph{exit boundary}
and, with $K(\cdot,\alpha)$ as basis for an $h$-transform, of the resulting transformed 
chain as the process conditioned on the value $\alpha\in\partial S$ 
for the limit $X_\infty$. Further, the relation to time reversal and the two-sidedness of
the Markov property outlined above shows that convergence $y_n\to\alpha\in\partial S$
in the Doob-Martin topology is equivalent to the convergence in distribution
of the initial segments of the conditioned chain, 
\begin{equation*}
     \cL[X_1,\ldots,X_k|X_n=y_n]\;\todistr\; \cL[X_1,\ldots,X_k|X_\infty=\alpha]
                         \quad\text{for all } k\in\bN.
\end{equation*}
Finally, we note that the topological construction depends on the chain via the 
associated cotransitions only. In the familiar forward case, the distribution of
$X$ is specified by the distribution of the first variable $X_1$ and the transition
function $p$. Here we may consider the distribution of $X_\infty$ as given,
and then use the cotransitions to obtain the distribution of $X$ as a mixture of
the distributions $\cL[X|X_\infty=\alpha]$, $\alpha\in\partial S$.

We now return to processes $(\Pi_n)_{n\in\bN}$ of permutations. This is a special case 
of the above framework, with $S=\bS$ and $S_n=\bS_n$. First, for $\PiF$
it is clear from the transition mechanism that, for $\sigma\in\bS_k$ and $\pi\in\bS_n$, 
the conditional probability  $P(\PiF_k=\sigma|\PiF_n=\pi)$ has value 1 
if the (order) restriction 
of $\pi$ to $[k]$ is equal to $\sigma$,
and that it is 0 otherwise. As functions on $\bS$ these are 
the same as the functions that generate the projective topology, see Section~\ref{sec:subs}.
Of course, convergence in the projective topology also follows directly 
from the fact that $\PiF$ is consistent in the sense that, for all $n\in\bN$, 
the restriction of $\PiF_{n+1}$ to $[n]$ is equal to $\PiF_n$.
At the other end of the spectrum, with the chain consisting 
of independent random variables, 
we obtain the one-point compactification, where each sequence of permutations that does
not have a limit point in $\bS$ itself converges to the single point $\infty$.
An intermediate case is $\PiH$ or, more generally, the copula-based 
process introduced in part (c) of Theorem~\ref{thm:Hoppen}. Recall that $\PiH$ appears
with the independence copula $C(x,y)=x\cdot y$.

\begin{theorem}\label{thm:DM}\label{DM} The process $\Pi=(\Pi_n)_{n\in\bN}$ generated 
by a copula~$C$ is a Markov chain with transition probabilities
\begin{equation}\label{eq:htransC}
   P(\Pi_{n+1}=\tau|\Pi_n=\sigma)\, 
           = \, \frac{(n+1)\,t(\tau,C)}{t(\sigma,C)} \, \pH(\sigma,\tau),
                          \quad n\in\bN,\, \sigma\in\bS_n,\tau\in\bS_{n+1},
\end{equation}
and the associated Doob-Martin topology is the same as the pattern frequency topology.
Further, $\Pi$ is an $h$-transform of $\PiH$.
\end{theorem}

\begin{RGproof}
We consider $\PiH$ first. By~\eqref{eq:cotransPiH} the functions in~\eqref{eq:funcDM} 
evaluate to $f_\sigma(\pi)=t(\sigma,\pi)$, which implies that the 
Doob-Martin topology coincides with the pattern frequency topology. 

For an arbitrary copula $C$
the permutations are generated by a sequence of independent random vectors $(X_i,Y_i)$, $i\in\bN$,
with distribution function $C$, and $\Pi_n$ is a deterministic function of the first $n$ pairs. 
The conditional distribution of $\Pi_n$ given all $\Pi_k$, $k>n$, thus depends only on $\Pi_{n+1}$
and the pair $(i,j)$ in the permutation plot of $\Pi_{n+1}$ that has been generated by 
$(X_{n+1},Y_{n+1})$. As the initial segments are invariant in distribution under permutations,
each $i\in [n+1]$ is equally likely. Taken together 
this shows that $(\Pi_n)_{n\in\bN}$ is a Markov chain with cotransitions 
that do not depend on $C$. This completes the proof of the topological assertion.

In order to obtain the expression for the one-step transition probabilities we note
that, by definition of $t(\cdot,C)$,
\begin{equation*}
    P(\Pi_k=\sigma) \,=\,  t(\sigma,C)\quad\text{for all }k\in\bN,\, \sigma\in\bS_k.
\end{equation*}
Using this together with~\eqref{eq:transPiH}, \eqref{eq:cotransPiH} and some elementary 
manipulations leads to~\eqref{eq:htransC}.

For the proof of the last assertion we need the extended Martin kernel $K(\cdot,\cdot)$
associated with $\PiH$.
First, by the familiar arguments,
\begin{equation*}
   K(\sigma,\pi)\, =\, \frac{P(\PiH_k=\sigma|\PiH_n=\pi)}{P(\PiH_k=\sigma)}   
                \; =\; k!\, t(\sigma,\pi)
\end{equation*}
for all $\sigma\in\bS_k,\pi\in\bS_n$ with $k\ge n$. Suppose now that $(\pi_n)_{n\in\bS}$ 
is such that $\pi_n\to C\in\partial\bS$. As the extension of $K(\sigma,\cdot)$ to $\bar \bS$ 
is continuous, it follows that 
\begin{equation*}
    K(\sigma,C) \,=\,\lim_{n\to\infty}K(\sigma,\pi_n) \,=\, k!\lim_{n\to\infty}t(\sigma,\pi_n) 
                \, =\, k!\; t(\sigma,C),
\end{equation*}
where we have used Theorem~\ref{thm:Hoppen}\,(a) in the last step. 
Using this~\eqref{eq:htransC} can be written as
\begin{equation*}
   P(\Pi_{n+1}=\tau|\Pi_n=\sigma)\, 
           = \, \frac{K(\tau,C)}{K(\sigma,C)} \, P\bigl(\PiH_{n+1}=\tau|\PiH_n=\sigma\bigr)
\end{equation*}
for all $n\in\bN$, $\sigma\in\bS_n$ and $\tau\in\bS_{n+1}$.
\qed
\end{RGproof}

We presented a purely topological approach to state space augmentation that 
is based on an embedding of $\bS$ into $[0,1]^\cF$ with a suitable family $\cF$ of functions $f:\bS\to [0,1]$ 
and the subsequent use of Tychonov's theorem. This is a variant of the famous  Stone-\v Cech
compactification; see e.g.~\cite[p.152]{Kelley}. Alternatively, we could endow $\bS$ with 
a metric $d$ and then  pass to the completion  of $(\bS,d)$. A suitable choice of $d$, based on  
the Martin kernel, leads to a totally bounded metric space that is topologically
homeomorphic to the compactification discussed above.

\section{Applications}\label{sec:app}
We give a range of examples where permutations and their limits are important,
from classical statistics to stochastic modeling to theoretical computer science.

\subsection{Tests for independence}\label{subsec:indep}
If the $X$- and $Y$-variables are independent then we obtain the copula 
$C(x,y)=x\cdot y$, $0\le x,y\le 1$. For the asymptotic frequencies this means that,
for all $k\in\bN$ and $\sigma\in\bS_k$,
\begin{align*}
   t(\sigma,C)\ 
      &=\  k!\, P(X_1<X_2<\ldots <X_k,\, Y_1<Y_2<\ldots<Y_k),               
\end{align*}
as a permutation of the $Y$-variables alone does not change the joint distribution of the
pairs $(X_1,Y_1),\ldots,(X_n,Y_n)$. In particular, all patterns of the same length have
the same asymptotic frequencies, so that $\, t(\sigma,C)=1/k!\, $ for all $\sigma\in\bS_k$.
 
This may be used for testing independence.
A famous classical procedure is based on  Kendall's $\tau$-statistic
which, in the notation introduced in Section~\ref{sec:subs}, may be written
as $\tau=\bigl(t(12,\Pi_n)-t(21,\Pi_n)\bigr)/ \binom{n}{2}$ and thus 
makes use of the two complementary patterns of length two. A central limit theorem 
for $U$-statistics leads to a test for independence with asymptotically
correct level, see e.g.~Example 12.5 in~\cite{vdV}. 
In our running data example a total of 61 of the 120 pairs are concurrent: Clearly, the test 
based on Kendall's $\tau$ would not reject the hypothesis of independence for these data.

An analogous procedure based on longer patterns requires a multivariate central 
limit theorem for vectors of $U$-statistics.
We consider patterns of length three, with a view towards
the values given in Table~\ref{tab:staedte}. For the asymptotic
covariances we use~Theorem~12.3 in \cite{vdV} and thus have to compute
$\zeta(\sigma,\tau):=3^2\,\cov(h_\sigma(Z_1,Z_2,Z_3),h_\tau(Z_1,Z'_2,Z'_3))$ 
for $\sigma,\tau\in\bS_3$, where 
$Z_1,Z_2,Z_3,Z'_2,Z'_3$ are independent and uniformly distributed on the unit square.   
For this we first determine the conditional distribution
$\cL[\Pi_3|X_1=x,Y_1=y]$ under the assumption of independence of horizontal and vertical 
positions. If $\pi=123$, for example, we either have  
$X_2<x,Y_2<y,X_3>x,Y_3>y$
or  $X_3<x,Y_3<y,X_2>x,Y_2>y$ if $(x,y)$ is `in the middle' of the rank plot, and similar
conditions apply in the other two cases. With $\phi_\sigma(x,y):= P(\Pi_3=\sigma|X_1=x,Y_1=y)$
elementary calculations along these lines lead to 
\begin{align*}
   \phi_{123}(x,y) \; &= \; 2x(1-x)y(1-y) + x^2 y^2/2 + (1-x)^2(1-y)^2/2,\\
   \phi_{132}(x,y) \; &= \; x^2y(1-y) + x(1-x)y^2 + (1-x)^2(1-y)^2/2,\\
   \phi_{213}(x,y) \; &= \; x(1-x)(1-y)^2 + (1-x)^2y(1-y) + x^2y^2/2,\\
   \phi_{231}(x,y) \; &= \; x(1-x)y^2 + (1-x)^2y(1-y) + x^2(1-y)^2/2,\\
   \phi_{312}(x,y) \; &= \; x(1-x)(1-y)^2 + x^2y(1-y) + (1-x)^2y^2/2,\\
   \phi_{321}(x,y) \; &= \; 2x(1-x)y(1-y) + x^2(1-y)^2/2 + (1-x)^2y^2/2.\\
\end{align*}
Note that $\phi_\sigma(x,y)=\phi_{\sigma^{-1}}(y,x)$ for all $\sigma\in\bS$, $0< x,y< 1$.
With these conditional probabilities we get 
\begin{align*}
 \cov&(h_\sigma(Z_1,Z_2,Z_3),h_\tau(Z_1,Z'_2,Z'_3))\\
      &=\ E\bigl(h_\sigma(Z_1,Z_2,Z_3)h_\tau(Z_1,Z_2',Z_3')\bigr)
                    \; -\;  \bigl(E h_\sigma(Z_1,Z_2,Z_3)\bigr)
                              \bigl(E h_\tau(Z_1,Z_2,Z_3)\bigr)\\
      &=\ \int_0^1\!\!\int_0^1 \phi_\sigma(x,y)\phi_\tau(x,y)\, dx\, dy    
                     \; -\; \Bigl(\int_0^1\!\!\int_0^1 \phi_\sigma(x,y)\, dx\, dy\Bigr)
                            \Bigl(\int_0^1\!\!\int_0^1 \phi_\tau(x,y)\, dx\, dy\Bigr).
\end{align*}
The observed value 0.228 for $\sigma=312$ deviates notably from the theoretical $1/6$. 
For this permutation the asymptotic variance is $\zeta=9(19/600-1/36)=7/200$. The 
corresponding asymptotic $p$-value for the observation would be the probability that
a standard normal variable exceeds the value
\begin{equation*}
     \sqrt{16} \, \sqrt{\frac{200}{7}}\, \biggl(0.228-\frac{1}{6}\biggr) 
                \ \approx \ 1.311,
\end{equation*} 
which is about $0.094$ and thus significant at the level $0.1$.
This, of course, is an act of `data snooping', but
it should be clear from the above how to obtain the $6\times 6$ asymptotic
covariance matrix for the vector $(t(\sigma,\Pi_n))_{\sigma\in\bS_3}$ of
pattern frequencies under the hypothesis of independence, and how then to apply
an appropriate test with the full data set in Table~\ref{tab:staedte}.
 
The use of ever longer patterns is a straightforward if computationally
demanding generalization. In this context another interesting connection between
discrete mathematics and statistics appears: A (deterministic) sequence 
$(\pi_n)_{n\in\bN}$ of permutations with $|\pi_n|\to\infty$ is said to be
\emph{quasirandom} if $t(\sigma,\pi_n)$ converges to $1/|\sigma|!$ for all 
$\sigma\in\bS$. The (random) sequence $(\Pi_n)_{n\in\bN}$ generated by the independence
copula has this property with probability one, so such sequences exist (there is an
obvious parallel with normal numbers). Is it enough for quasirandomness if the limit
statement holds for all permutations of some fixed length $k$? 
It follows from \eqref{eq:cotranspattern} that this `forcing property'  would 
then automatically hold for all $j\ge k$ and, clearly, $k=2$ is too small.
This question has been answered in~\cite{Kral}, who showed that quasi\-ramdomness follows 
with $k=4$, but not with $k=3$. Their proof of the sufficiency part is based on the
permuton approach discussed in Section~\ref{sec:subs}. 
In fact, as every limit point $C$ must satisfy $t(\sigma,C)=1/4!$ 
for all $\sigma\in\bS_4$, the statement boils down to a characterization property 
of the independence copula. Remarkably, in a statistical context, such a 
characterization had been obtained much earlier in~\cite{Yana}, and it leads
to an extension of Kendall's $\tau$ that is consistent against
all alternatives. In both areas, related questions are active fields of research.
In connection with quasirandomness we refer to~\cite{sixperm},
for the statistical side see e.g.~\cite{Drton}.

\subsection{Delay models}\label{subsec:delay} 

Objects arrive at a system at times $U_1, U_2, \ldots$ and leave at times 
$U_1 + X_1, U_2 + X_2, \ldots$, where we assume that the arrivals are independent and 
uniformly distributed on the unit interval, that the delay times are independent with
distribution function $G$, and that arrival and delay times are independent. Let 
$\Pi_n$ be the random permutation that connects the orders of the first $n$ arrivals 
and departures. For example, in a queuing context, arriving customers might receive 
consecutively numbered tickets which they return on departure, and the tickets are put 
on a stack.

In~\cite{BaGrDeldat} the resulting permutons, the delay copulas, are
discussed. From a probabilistic and statistical point of view the information about
the delay distribution contained in the permutations is of interest, together with
the question of how to estimate (aspects of) the delay distribution. Another point of interest 
are second order approximations or, in general terms, the transitions from a strong
law of large number to a central limit theorem. For the set of patterns of a fixed
size this has already been important in Section~\ref{subsec:indep}. 
In~\cite{BaGrDeldat} this is taken further by considering 
$\sqrt{n}(\Pi_n(\sigma)-t(\sigma,C))_{\sigma\in \bS}$ as a stochastic process 
with $\sigma\in\bS$ as `time parameter' and then establishing a functional 
central limit theorem.

\subsection{Sparseness}\label{subsec:sparse} 
The delay models considered in Section~\ref{subsec:delay} can be related to the
$M/G/\infty$ queue where customers arrive at the times of a Poisson process
with intensity $\lambda$. There are infinitely many servers so, strictly speaking, 
there is no queuing. 
The limits refer to an increasing arrival rate and thus generate a 'dense' sequence
of permutations as the service time distribution is kept fixed.  
Here we regard the standard $M/G/1$ model with  only one server and
fixed arrival rate; for convenience we start with an empty queue at time 0. From 
the first $n\in\bN$ arrivals and departures  we still obtain a sequence 
$(\Pi_n)_{n\in\bN}$ of random permutations, with $\Pi_n\in\bS_n$ for all $n\in\bN$.
We argue that this sequence is `sparse', and that the pattern frequency topology 
does not provide much insight.

We assume that the service time distribution has finite first moment $\mu$ and that
$\rho <1$, where $\rho:=\lambda\mu$ is known as the \emph{traffic intensity}. The queue is then
stable and the queue length process $Q=(Q_t)_{t\ge 0}$, with $Q_t$ the number of customers
in the system at time $t$, can be decomposed into alternating busy and idle periods.
These are independent, and the individual sequences of busy and idle pieces of $Q$ 
are identically distributed. Obviously,  arrival and departure times of any particular 
customer belong to the same busy period. Let 
$K_i$ be the number of customers in the $i$th period. We recall the 
definition of the direct sum of finite permutations: 
For $\sigma\in\bS_k$, $\tau\in\bS_l$, we obtain the sum $\sigma\oplus \tau\in\bS_{l+k}$
in one-line notation as
\begin{equation*}
        \sigma\oplus \tau =(\sigma_1,\sigma_2,\ldots,\sigma_k,\tau_1+k,\tau_2+k,\ldots,\tau_l+k).   
\end{equation*} 
Suppose now that $K_1=k$. Then, for all $n\ge k$, $\Pi_n(i)\le k$ if $i\le k$, and
$\Pi_n(i)>k$ if $k<i\le n$, hence $\Pi_n=\sigma\oplus \tau$ with some $\sigma\in\bS_k$. 
Iterating this argument we see that, as $n\to\infty$, $\Pi_n$ decomposes into 
subpermutations related to the sequence of busy cycles. Also, for any given $i$, 
$\Pi_n(i)$ remains constant as soon as $n$ exceeds the number
of the first customer no longer in the same busy period as customer $i$. 
Taken together this leads to the following result.

\begin{proposition}\label{prop:queue1}
Let $(\Pi_n)_{n\in\bN}$ be the sequence of permutations generated by an $M/G/1$ queue 
with traffic intensity $\rho<1$. Then $\ \Pi_n\to \Pi_\infty$ almost surely in the 
projective topology, with 
$\;\Pi_\infty\eqdistr\bigoplus_{i=1}^\infty \Psi_i$, where
$\Psi_i\in\bS$, $i\in\bN$, are independent and identically distributed 
and $|\Psi_i|=K_i$.
\end{proposition}

The distribution $\cL(\Pi_\infty)$ of the limit $\Pi_\infty$ is fully characterized by
$\cL(\Psi_1)$, which may in turn be specified by $\cL(K_1)$ and the conditional distribution 
$\cL(\Psi_1|K_1)$. Nothing is known (to me) about the latter, apart from the obvious fact
that it is concentrated on the subset of $\bS_{K_1}$ of permutations that cannot be written 
as a direct sum. 

The sequence $(\Pi_n)_{n\in\bN}$ in Proposition~\ref{prop:queue1} is contained in 
the compactification of $\bS$ that we used in the `dense' situation. 
Under a mild moment condition we obtain the corresponding set of limit points. 

\begin{proposition}\label{prop:queue2}
Suppose that, in addition to the assumptions in Proposition~\ref{prop:queue1}, 
the service time distribution has finite third moment. Then,
in the pattern frequency topology, $\Pi_n\to C_\infty$ almost surely,
where the permuton $C_\infty$ is given by $C_\infty(u,v)=\min\{u,v\}$ for all 
$u,v\in[0,1]$.
\end{proposition}

\begin{RGproof} 
The path-wise argument used for Proposition~\ref{prop:queue1} 
shows that for any pair $(i,j)$ with $i<j$ and 
$\Pi_n(i)>\Pi_n(j)$ customers $i$ and $j$ must belong to the same busy period. In particular,
\begin{equation*}
   \#\bigl\{1\le i<j\le n:\, \Pi_n(i)>\Pi_n(j)\bigr\}   
                  \, \le \,  n \, M_n,\quad\text{with } M_n:= \max\{K_1,\ldots,K_n\},
\end{equation*}
so that the relative frequency of the inversion $\tau=21$ satisfies 
$t(\tau,\Pi_n) \le 2\,M_n/(n-1)$. It follows from Section 5.6, p158 in \cite{CoxSmith}
that $EK_1^3<\infty$ under the above moment assumption. 
By a standard argument this implies $M_n/n\to 0$ almost surely.

Now let $\sigma$ be a permutation that contains an inversion. The relative frequency
$t(\sigma,\pi)$ of $\sigma$ in $\pi\in\bS_n$, $n\ge k:=|\sigma|$, can be interpreted as the
probability that $\pi[A]=\sigma$, if $A\subset [n]$ with $\# A=k$ is chosen uniformly at random;
see also the remarks at the beginning of Section~\ref{sec:subs}.
Further, choosing $\{i,j\}\subset [k]$ leads to an inversion with probability $p>0$, where $p$ 
does not depend on $n$. In the two-stage experiment, with independent steps, we would thus 
find an inversion in $\Pi_n$ with probability $p\cdot t(\sigma,\Pi_n)$. This is bounded from 
above by $t(\sigma,\Pi_n)$, and the first step of the proof therefore implies
that $t(\sigma,C) =\lim_{n\to\infty} t(\sigma,\Pi_n)=0$ for all permutations that contain 
an inversion. As the subsampling must lead to \emph{some} element of $\bS_k$ we see that
any limit permuton $C$ must assign the value~1 to each identity permutation. This holds for 
the copula in the assertion of the theorem, and an appeal to uniqueness of the limit 
completes the proof.
\qed 
\end{RGproof}

Thus, from the pattern frequency point of view, there is asymptotically no difference 
between $M/G/1$ and the queue where customers depart in order of their arrival, such as
in a system with immediate service and constant service times.

\subsection{Pattern avoiding and stack sorting}\label{subsec:avoid}
In~\cite{Bass} the authors deal with \emph{separable permutations}, which are those
that avoid the two patterns $2413$ and $3142$. This class is of interest as it appears
in connection with stack sorting, for example. (The permutation in our data example is not separable; see the red dots in Figure~\ref{fig:Staedte}.)

Starting with the sets $\AV_n:= \{\pi\in\bS_n:\, t(2413,\pi)=0,\, t(3142,\pi)=0\}$ the main 
question is the limit distribution (in the weak topology related to the pattern frequency
topology) of $\Pi_n\sim \unif(\AV_n)$ as $n\to\infty$. Whereas we often
obtain a limit that is concentrated at one point, such as the independence copula
in connection with the unrestricted case where $\Pi_n=\unif(\bS_n)$ in 
Section~\ref{subsec:indep}, or in the delay context in Section~\ref{subsec:delay},
it turns out that the limit distribution for 
uniformly distributed separable permutations is `truly random'; in fact, this distribution
is an interesting object on its own.
The analysis in~\cite{Bass} is based on the fact that separable permutations may be coded 
by a specific set of trees, where the operation $\oplus$ mentioned in 
Section~\ref{subsec:sparse} together with its `skew' variant $\ominus$ plays a key role. 
I do not know at present if this
fits into the Markovian framework. Is there a Markov chain adapted to $\bS$ 
with marginal distributions uniform on $\AV_n$ and such that the Doob-Martin approach 
leads to the pattern frequency topology?

Again, this is an active area of research. We refer the reader to \cite{JansonAlg} for 
a recent contribution, which is of particular interest in the context of the present 
paper because of its use of a variety of probabilistic techniques.

\section{Complements}\label{sec:compl}
We sketch some related developments, restricting references essentially to 
those that can serve as entry points for a more detailed study.

\subsection{Other combinatorial families}\label{subsec:CombFam}
As pointed out in~\cite{Hopp} the construction of permutons as limits of permutations via
pattern frequencies has been influenced by the earlier construction of graphons
as limits of graphs via subgraph frequencies; see~\cite{Lov} for a definite account
and, in addition, \cite{DiaconisJanson} for a probabilistic view.
It was observed in~\cite{GrDMTCS} that, with a suitable Markov chain, the topology based
on subgraph frequencies is the same as the Doob-Martin topology arising from the 
associated boundary, in analogy to Theorem~\ref{thm:DM} above. A notable difference
between the permutation and the graph situation is that graphs are considered modulo
isomorphisms. In~\cite{GrDMTCS} a randomization step was used to take care of this,
\cite{Hage} showed that one can work directly with the equivalence classes.

Markov chain boundary theory has been used in connection with limits of random discrete 
structures in a variety of combinatorial families other than $\bS$; 
see~\cite{GrSemBerKMK} for an elementary introduction and references. 
We single out the case of binary trees, where $\bB_n$ is 
the set of binary trees with $n$ (internal) nodes and $\bB=\bigsqcup_{n=1}^\infty \bB_n$ is 
the graded state space.
Processes $X=(X_n)_{n\in\bN}$ with values in $\bB$ arise naturally in 
two different situations.
First, if we apply the BST (binary search tree) algorithm sequentially to a sequence of
independent and identically distributed real random variables with continuous distribution 
function, and secondly, in the context of R\'emy's algorithm, where $X_n$ is uniformly 
distributed on $\bB_n$. Both processes are Markov chains adapted to $\bB$, but 
the associated transition graphs and, consequently, the cotransitions and the Doob-Martin compactifications of $\bB$, are quite different. For search trees the role 
of pattern or subgraph frequencies is taken over by relative subtree sizes, and the 
permuton or graphon analogues are the probability distributions on the set of infinite 0-1 
sequences; see~\cite{EGW1}. In the R\'emy case a class of real trees appears, together
with a specific sampling mechanism; see~\cite{EGW2}.

\begin{figure}
\setlength{\abovecaptionskip}{-.25cm}

\includegraphics[scale=.5]{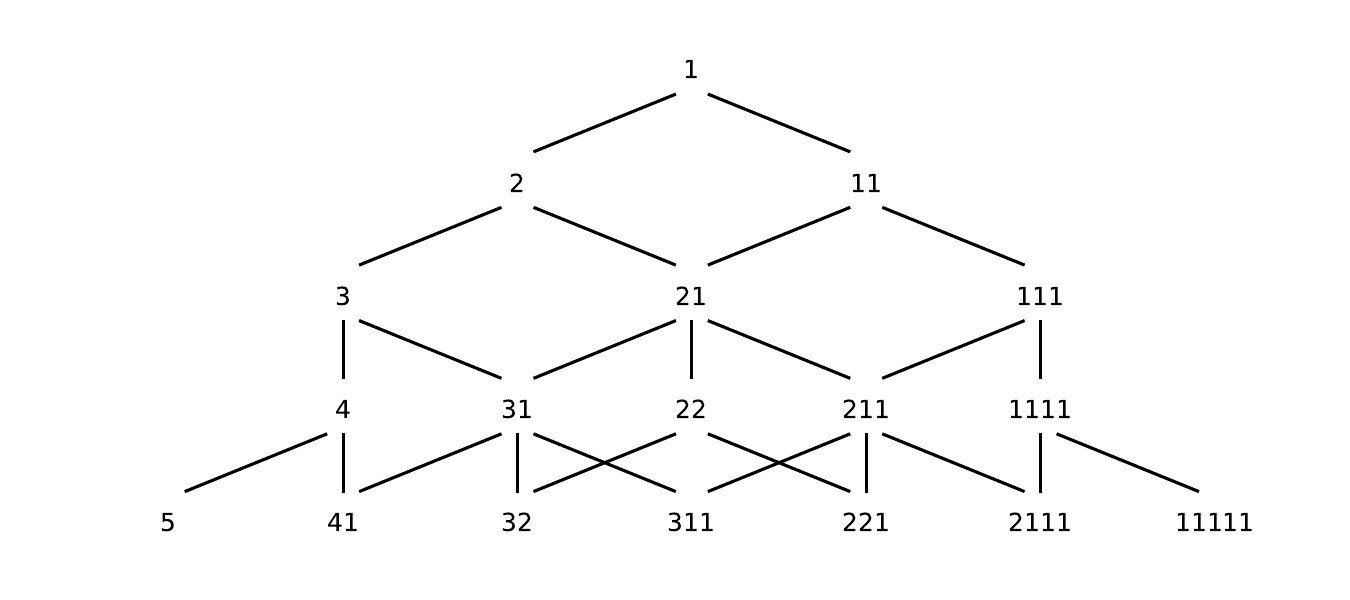}
\vspace{-.1cm}
\caption{The transition graph for integer partitions}
\label{fig:Kingman}
\end{figure}

Another combinatorial family that is related in several ways to the topics discussed 
in the present paper 
is the family $\bY=\bigsqcup_{n\in\bN} \bY_n$  of (number) partitions. 
Here $\lambda=(\lambda_1,\ldots,\lambda_k)\in\bY_n$ is a \emph{partition}
of $n$ if $\lambda_i\in\bN$,
$\lambda_1\ge \cdots\ge \lambda_k$, and $\sum_{i=1}^k \lambda_i=n$. We abbreviate this to
$\lambda\vdash n$. 
Each $\sigma\in\bS_n$ defines a partition $\lambda\vdash n$ via the length of its cycles, 
in decreasing order. The right hand side of~\eqref{eq:Foata}, for example, leads to
$\lambda=(5,3,1)$.
On $\bY$ a partial order can be defined as follows: If 
$\lambda=(\lambda_1,\ldots,\lambda_k)\in\bY_m$ and $(\eta_1,\ldots,\eta_l)\in\bY_n$
with $m\le n$ then $\lambda\le \eta$ means that $k\le l$ and $\lambda_j\le \eta_j$ for
$1\le j\le k$. Figure~\ref{fig:Kingman} shows the first five levels of the associated
Hasse diagram. We obtain weights for the edges of the graph by `atom removal': For example,
to go from $(3,2)$ to the two predecessors $(3,1)$ and $(2,2)$ respectively there are two 
possibilities to decrease $2$ to $1$ in the first and three to decrease $3$ to~$2$ in 
the second case. Normalizing these values 
so that they sum to~1 we obtain the cotransitions for a family of Markov chains 
adapted to $\bY$.  

Suppose now that $(\lambda(n))_{n\in\bN}$ is a sequence in $\bY$ 
where we again assume for simplicity
that $\lambda(n)\in\bY_n$ for all $n\in\bN$. It turns out that the Doob-Martin topology
associated with these cotransitions 
is equivalent to the convergence of the normalized partition parts $\lambda(n)_i/n$ as 
$n\to \infty$, for all $i\in\bN$, and that 
$\partial \bY=:\partial \bY_{\text{\tiny\rm c}}$ is homeomorphic to the
space of weakly decreasing sequences $(\alpha_i)_{i\in\bN}$ with 
$\sum_{i=1}^\infty\alpha_i\le 1$,
endowed with the trace of the product topology on $[0,1]^\infty$.

\subsection{The cycle view}\label{subsec:cycle}
Above we worked with the order isomorphism aspect of permutations
and the associated partial order on $\bS$ given by pattern containment. 
Using cycles instead leads 
to another partial order on $\bS$ where, for $\sigma\in\bS_m$ and $\tau\in\bS_n$
with $m\le n$, the relation $\sigma\le \tau$ means that in the cycle notation, $\sigma$ arises 
from $\tau$ by deleting all numbers greater than $m$ (some rearrangement may be necessary 
to obtain the standard notation). As an example we consider the two 
sides of~\eqref{eq:Foata} again: The order restriction to the set 
$[4]$ in the one-line notation 
on the left leads to $3241\in\bS_4$, which is $(134)(2)$ in standard cycle notation. The
cycle restriction to $[4]$ of the permutation on the right hand side is $(1)(432)$ in 
standard notation, which is $1423$ in one-line notation. In particular, Foata's correspondence 
is not consistent with the two notions of restriction.

Nevertheless, the correspondence, together with the construction of $\PiF$ on the left hand 
side, can be used to motivate a Markov chain $\PiC=(\PiC_n)_{n\in\bN}$ that is adapted to 
$\bS$ and compatible with this order, and has $\PiC_n\sim\unif(\bS_n)$ for all $n\in\bN$.
In this construction, known as the \emph{Chinese restaurant process}, cycles are interpreted 
as circular tables and customer $n+1$ selects one of the $n$ already 
seated customers as right neighbor (successor in the cycle notation),  or starts a new table, 
where the $n+1$ possibilities are chosen with the same probability. 
Note that $\PiC$ is projective in the sense that previous values are functions of the
present state. 

The cycle view is closely related to number partitions: Mapping a permutation to its
ordered cycle lengths leads to a function $\Psi:\bS\to\bY$ that preserves the grading and the
respective partial orders. 
Let $\LaC=(\LaC_n)_{n\in\bN}$ be defined by $\LaC_n:=\Psi(\PiC_n)$ for all $n\in\bN$.
Then $\LaC$ is again a Markov chain, and it has the cotransitions introduced in
Section~\ref{subsec:CombFam} for the family $\bY$. Hence, almost surely as $n\to\infty$,
$\LaC_n$ converges  to an element $(\alpha_i)_{i\in\bN}$
of $\partial \bY_{\text{\tiny \rm c}}$ in the sense that 
$\lim_{n\to\infty}(\LaC_n)_i/n=\alpha_i$ for all $i\in\bN$.
The distribution of the limit, which is not concentrated at one point of the boundary, 
can be obtained by \emph{stick breaking}: Starting with
independent and $\unif(0,1)$-distributed random variables $U_i$, $i\in\bN$, 
let $(V_i)_{i\in\bN}$ be
defined by  $V_1=U_1$ and $V_{n+1}=(1-V_n)U_{n+1}$ for all $n\in\bN$.
Putting these in decreasing order provides a random element of 
$\partial\bY_{\text{\tiny\rm c}}$ that has the same distribution as $\LaC_\infty$.  

Note that this approach refers to the big cycles. Locally, at the small end, the limit
can be obtained from~\eqref{eq:cycrepr}: For $i\in\bN$ fixed, the counts of cycles of 
length $k$, $1\le k\le i$, are asymptotically independent and Poisson distributed with
parameters $1/k$.

The cycle notation displays a permutation as a product of cyclic permutations
that act on disjoint parts of the base set. In contrast to the general case such products 
are commutative. The cycle view is closely connected to the group aspects of permutations, 
in particular to the representation theory of finite groups.
Many probabilists and statisticians became aware of this connection 
and the resulting `non-commutative Fourier analysis' through~\cite{Diaconis}; more recent and
extensive presentations are~\cite{CSHarm} and~\cite{Mel}. We briefly summarize some
aspects that are relevant to the asymptotics of large permutations. 
All representations below refer to base field~$\bC$, $\GL(V)$ denotes the general linear 
(or automorphism) group of the complex vector space~$V$, and $U(n)$ is the group of 
unitary $n\times n$-matrices. 

As with compactifications, we begin with an embedding:
We regard the elements $g$ of a finite group $G$ as elements of the vector space 
$\bC^G$ of functions $f:G\to \bC$ via $\sigma\mapsto\delta_\sigma$, 
where $\delta_\sigma(\sigma)=1$ and $\delta_\sigma(\tau)=0$
if $\tau\not=\sigma$.  With the composition of $\sigma,\tau\in G$ written as $\sigma\tau$ 
we then define convolution on $\bC^G$ by
\begin{equation*}
   f\star g (\sigma) = \sum_{\tau\in G} f(\sigma \tau^{-1})\, g(\tau), 
                                      \quad \sigma\in G.
\end{equation*}
Seen as multiplication, convolution  makes $\bC^G$ an algebra, the \emph{group algebra}
$\bC[G]$ associated with~$G$. Also, $(\sigma,f)\mapsto \sigma.\, f:=\delta_\sigma \star f$
defines an action of $G$ on $V$, and 
the function $\rho:G\to \GL(V)$, $\rho(\sigma)(f):= \sigma .\, f$, is a 
group homomorphism. With respect to the canonical basis \hbox{$\{\delta_\sigma:\,\sigma\in G\}$}
the values $\rho(\sigma)$ are all represented by unitary matrices. 
Thus, somewhat reminiscent of the fact that
every finite group $G$ with $n$ elements is a subgroup of $\bS_n$, every such $G$ is a
subgroup of $U(n)$, both up to group isomorphism. This is known as the (left)
\emph{regular representation} of $G$. 
More generally, a representation $(\rho,W)$ of $G$ consists of a vector space $W$ 
and a homomorphism $\rho$ from $G$ to $\GL(W)$. 
For example, representing permutations
by their permutation matrices as in~\eqref{eq:permmatrix}, we obtain the 
\emph{permutation representation} of $G=\bS_n$.

In a representation $(\rho,W)$ of $G$ a subspace $U$ of $W$ is invariant 
if $\rho(\sigma)(U)=U$ for all $\sigma\in G$, and the representation is 
\emph{irreducible} if only $U=\{0\}$ and $U=W$ are invariant. Necessarily, 
the range of $\rho$ is then the full group $\GL(U)$. It is a crucial fact that
the regular representation $(\rho,\bC[G])$ can be decomposed into the direct sum
$(\bigoplus_{i\in I} \rho_i,\bigoplus_{i\in I} W_i)$ of irreducible representations 
$(\rho_i,W_i)$. Let $d_i:=\dim W_i$. For later use we note that, with $n$ the size of $G$,
counting dimensions leads to $\sum_{i\in I} d_i^2=n$.  
The \emph{character} $\chi:G\to\bC$ of a representation $\rho$ is given 
by $\chi(\sigma):=\Tr(\rho(\sigma))$, where $\Tr$ denotes the trace (which does not depend 
on the basis chosen to represent $\rho(\sigma)$ by a matrix). For the permutation 
representation $\chi(\sigma)$ is easily seen to be the number of fixed points of~$\sigma$.

We recall that general groups can be decomposed into conjugacy classes, which are 
the equivalence classes obtained when $\sigma,\tau$ are equivalent if $\sigma\pi=\pi\tau$
for some $\pi\in G$. A function  $f\in \bC[G]$
is a \emph{class function} or \emph{central} if it is constant on conjugacy classes. 
The characters of representations are such class functions, 
and it turns out that the characters of the irreducible representations constitute
a basis for the space (in fact, a convolution algebra) of class functions.  
Moreover, different irreducible representations lead to different characters. 
For $G=\bS_n$ the conjugacy classes correspond to cycle structures and may thus be represented
by partitions of $n$, i.e.\ elements of $\bY_n$. 
Remarkably, the characters of irreducible representations can also be parametrized by such 
partitions. 

This finishes our excursion, and we return to permutation asymptotics. The algebraic 
point of view mainly refers to characters and thus to partitions, where it leads to
a close relative of the Markov chain $\LaC=(\LaC_n)_{n\in\bN}$ discussed above.
For this, our main source is the book \cite{Kerov}; see also the references given there.

It is convenient to augment the 
permutations of $[n]$ by fixed points so that a bijection of $\bN$ is obtained that
leaves all $k>n$ invariant. We continue to write $\bS_n$ for these `padded' versions,
and regard them as an increasing sequence of 
subgroups of the group $\bS_\infty$ that consists of all finite
bijections $\sigma:\bN\to\bN$, meaning that $\#\{k\in\bN:\, \sigma(k)\not=k\}<\infty$.
A representation $\rho_{n+1}$ of $\bS_{n+1}$ then leads to a 
representation $\rho_n$ of $\bS_n$ by restriction. This will in general destroy 
irreducibility, but the character of $\rho_n$ is a class function and may thus be written 
as a linear combination of the irreducible characters on level $n$. It turns out that
the non-zero coefficients are all equal to one. Consider now the graph 
on $\bY$ with edges connecting the character of an irreducible representations on level 
$n+1$ to the contributing characters of the irreducible representations on level $n$.
This is the same graph as the transition graph obtained above for number partitions,
see Figure~\ref{fig:Kingman}, but the edge weights are now all equal to 1. The resulting 
cotransitions then amount to the rule that, from a current state $\eta\in\bY_n$, a state $\lambda\in\bY_{n-1}$ is chosen uniformly at random from its predecessors. 
The dimension numbers $d_\lambda$, $\lambda\in\bY$, turn out to be equal to
the number of paths in the diagram that lead from its root to $\lambda$. As a consequence
the cotransition associated with an edge $\{\lambda,\eta\}$, $\lambda<\eta$, of the 
transition graph is given by 
$d_\lambda/d_\eta$.

When does a sequence of characters $(\chi_n)_{n\in\bN}$ of irreducible representations or,
equivalently, a sequence of partitions $(\lambda(n))_{n\in\bN}$, converge in the
Doob-Martin topology associated with these cotransitions?
Or, in other words, what is the cycle view analog of the pattern frequency topology?
The answer requires one more definition: The \emph{conjugate} (or transpose) of a partition 
$\lambda=(\lambda_1,\ldots,\lambda_k)$ is given by $\lambda^\star=(\lambda^\star_1,\ldots,\lambda^\star_j)$ with $j=\lambda_1$ and
$\lambda^\star_l=\#\{m\in[k]:\, \lambda_m\ge l\}$. Convergence then means that the normalized
partition components of $\lambda(n)\in\bY_n$ \emph{and} $\lambda^\star(n)$ converge as 
elements of the unit interval: $\lambda(n)_i/n\to \alpha_i$, $\lambda^\star(n)_i/n\to\beta_i$ 
as $n\to\infty$, for all $i\in\bN$. The boundary is known as the \emph{Thoma simplex},
\begin{equation*}
 \partial\YT \, =\, \Bigl\{(\alpha,\beta)\in[0,1]^\bN\times [0,1]^\bN:\, 
                            \alpha_i\ge \alpha_{i+1}, \beta_i\ge \beta_{i+1}
                                      \text{ for all }i\in\bN, \ 
                          \sum_{i=1}^\infty (\alpha_i+\beta_i)\le 1 \Bigr\},
\end{equation*}
together with the topology of pointwise convergence. 
The boundary $\partial\bY_{\text{\tiny\rm c}}$ for the partition chain 
obtained above in connection with the Chinese restaurant process
can be identified with the compact subset of $\partial\YT$ that arises
if all $\beta$-values are equal to  zero.

A specific combinatorial Markov chain adapted to $\bY$ is given by the 
\emph{Plancherel growth process} $\LaP=(\LaP_n)_{n\in\bN}$, with transition probabilities 
\begin{equation*}
       P(\LaP_{n+1}=\eta|\LaP_n=\lambda)\, =\, \frac{d_\eta}{(n+1)d_\lambda}, 
           \quad \text{for all }\lambda\in\bY_n,\,\eta\in\bY_{n+1} \text{ with }\lambda<\eta,
\end{equation*}
and marginal distributions $P(\LaP_n=\lambda)= d_\lambda^2/n!$, $n\in\bN$,
$\lambda\in\bY_n$. It follows from the above dimension counting that the latter are indeed 
probability mass functions on $\bY_n$, and the calculation between Lemma~\ref{lem:Hoppen}
and Theorem~\ref{thm:harm}
confirms that $\LaP$ has the cotransitions associated with character restrictions. In contrast
to the partition chain $(\LaC_n)_{n\in\bN}$ arising in the cycle length context, 
the limit distribution for the Plancherel chain $(\LaP_n)_{n\in\bN}$ is concentrated 
at the one point of the boundary $\partial\YT$ given by $\alpha_i=\beta_i=0$ for all $i\in\bN$. 
The analysis of this chain leads to some quite
spectacular results, such as the arcsine law for the limit shape of 
Young diagrams, and the solution of Ulam's problem on the longest 
increasing subsequences in $\Pi_n\sim\unif(\bS_n)$.  
We refer to~\cite{Kerov} again, for Ulam's problem see also~\cite{Romik}. 

In connection with trees we had two different graph structures on $\bB$. Here the graphs
on $\bY$ are the same, but the weights are different. A generalization that includes both
the cycle model and the Plancherel process is given 
in~\cite{KOO}. There, an important aspect is the appearance of classical families of symmetric 
polynomials in the context of the respective extended Martin kernel.

\subsection{Exchangeability}\label{subsec:deFinetti}
In Section~\ref{sec:dyn} we worked with a general Markov chain $X=(X_n)_{n\in\bN}$
that is adapted to a graded state space $S=\bigsqcup_{n\in\bN} S_n$ and we used the 
cotransitions of $X$ to obtain a compactification $\bar S= S\sqcup \partial S$ 
of the state space, together with a limit variable $X_\infty$ for the $X_n$'s as $n\to\infty$. 
Each boundary value $\alpha\in\partial S$ induces a probability measure $P^\alpha$, which is 
the model for $X$ conditioned on $X_\infty=\alpha$.
Further, an arbitrary distribution $\mu$ on the Borel
subsets of $\partial S$ may be used to construct a Markov chain that has the
same cotransitions as $X$ and limit distribution 
$\cL(X_\infty)= \int P^\alpha\, \mu(d\alpha)$, the $\mu$-mixture of the $P^\alpha$'s. 
We note that this is the time reversal version
of the familiar `forward' situation where the distribution of a Markov chain is specified
by the distribution of the first variable and the (forward) transition probabilities. 
Here, however,
the distribution of $X_\infty$ is usually defined on a non-discrete space, and
there is of course no step from~$\infty$ to~`$\infty-1$'. 

The above shows some similarities to the archetypical exchangeability result, 
de Finetti's theorem:
We start with a sequence $X=(X_n)_{n\in\bN}$ of real random variables that is exchangeable 
in the sense that the distribution $\cL(X)$ of $X$ is invariant under all 
$\sigma\in\bS_\infty$, the set of finite 
permutations of $\bN$ introduced in Section~\ref{subsec:cycle}. We then
obtain a limit $M_\infty$ for the empirical distributions
$M_n:=\frac{1}{n}\sum_{i=1}^n\delta_{X_i}$ and it holds that
$\cL(X|M_\infty=\mu)=\cL(\tilde X)$ where $\tilde X=(\tilde X_i)_{i\in\bN}$ 
with $\tilde X_i$, $i\in\bN$, independent random variables with distribution $\mu$.
This area has developed into a major branch of modern probability, dealing with 
structural results for distribution families with certain invariance properties. 
\cite{KallSym} is a standard reference, but see also~\cite{AldousSF} and~\cite{Austin}.
An excellent introduction is also given in the first part of the lecture notes \cite{AustinLN}.

Exchangeability may provide an approach to the asymptotics of 
random discrete structures; see e.g.\ Section 11.3.3 in~\cite{Lov} for graphs 
and~\cite{EGW2} for R\'emy trees. (In the other direction, a proof for the classical 
de Finetti theorem can be given using Markov chain boundary theory; see~\cite{GGH}.) 
In their treatment of random permutations and 
permutons \cite{Hopp} emphasize the connection to random graphs and graphons. Our
aim here is to extend this to the exchangeability aspect.

In the subgraph frequency topology for sequences of (random) graphs 
the limits are described by a (possibly random)
\emph{graphon}, a measurable and symmetric function 
$W:[0,1]\times [0,1]\to [0,1]$. Given $W$, we use 
a sequence $U=(U_i)_{i\in\bN}$ of independent uniforms to construct a  random binary
$\bN\times\bN$-matrix $Z$ by choosing $Z_{ij}=Z_{ij}=1$ with probability $W(U_i,U_j)$, 
$1\le i < j<\infty$, independently for different pairs, and $Z_{ii}=0$ for all $i\in\bN$. 
This matrix is jointly exchangeable in the sense that, for all $\sigma\in\bS_\infty$,
$Z^\sigma:=(Z_{\sigma(i),\sigma(j)})_{i,j\in\bN}$
has the same distribution as $Z$. The upper $n\times n$-corner of $Z$ is the 
incidence matrix of a random graph $G_n$.
As $n\to \infty$, these converge in the subgraph frequency topology, and the limit is 
represented by $W$. The graphon may be regarded as an analogue of the
measure $M_\infty$ in the classical case of exchangeable sequences, and the second step 
as an analogue of sampling from $M_\infty$.

Is there a similar representation for random permutations, specifically for the
copula-based models in part (c) of Theorem~\ref{thm:Hoppen}? We need an infinite matrix
that represents the complete sequence, with the upper left $n\times n$-corner for 
the result of the first $n$ steps.  
Starting with independent random vectors $(X_i,Y_i)$, $i\in\bN$, with distribution 
function~$C$, we define an infinite 
random binary array $Z=(Z_{ij})_{i,j\in\bN}$ by
\begin{equation*}
     Z_{ij}=\begin{cases}
                 1, &\text{if }i<j\text{ and } Y_i<Y_j,\\
                 1, &\text{if }i>j\text{ and } X_i<X_j,\\
                 0, &\text{else.}  
            \end{cases}
\end{equation*}
For example, if the cities in our data set arrive in order of decreasing population
size (Berlin, Hamburg, Munich, Cologne $\ldots$) then the upper left $4\times 4$-corner
of $Z$ is
\begin{equation*}  
(Z_{ij})_{i,j=1}^4\; = \; \begin{pmatrix} 
           0&1&0&0 \\0&0&0&0 \\0&1&0&1 \\0&0&0&0 \\
\end{pmatrix} .
\end{equation*}
Within this group, there is one city 
with more citizens than Munich that is further to the east, and none of them is further to the
south.
The permutation $\Pi_n$ can be obtained from $(Z_{ij})_{1\le i,j\le n}$ as follows: With
\begin{equation*}
    \#\bigl\{1\le i\le n:\, X_j\le X_i\bigr\}\; 
                  =\ \sum_{i=1}^{j-1}Z_{ij} \,+\, (n-j-1) \;- \sum_{i=j+1}^n Z_{ij},
                            \quad 1\le j\le n,
\end{equation*}
we get the absolute ranks of the variables $X_i$, $1\le i\le n$, within $\{X_1,\ldots,X_n\}$,
and a similar formula holds for 
the $Y$-components. Then $\Pi_n$ can be determined from the two rank vectors as
in~\eqref{eq:rank}. Equivalently, we may work with an array $\Tilde Z$ indexed by subsets 
$\{i,j\}$ of $\bN$ of size two, $i<j$, by using the values $0,1,2,3$ for indicating the
relative position of $(X_j,Y_j)$ with respect to $(X_i,Y_i)$ through the four 
quadrants North-East (NE), SE, SW, and~NW.

If $\sigma:\bN\to\bN$ is strictly increasing, then the matrix $Z^\sigma$ with entries
$Z_{\sigma(i),\sigma(j)}$ depends on $(X_{\sigma(i)},Y_{\sigma(i)})$, $i\in\bN$, 
in the same (deterministic) way as $Z$ depends 
on $(X_i,Y_i)$, $i\in\bN$. As the subsequence $(X^\sigma,Y^\sigma)$ is equal in 
distribution to the  original sequence $(X,Y)$, the arrays $Z^\sigma$ and $Z$ are 
then equal in distribution, which shows that $Z$ is contractable. An analogous statement 
holds for the corresponding array $\Tilde Z$ indexed by size two subsets of $\bN$. 
It follows from the general Aldous-Hoover-Kallenberg
representation theorem that a contractable array can be extended to an exchangeable array; 
see Corollary 7.16 in~\cite{KallSym}.

In order to construct a suitable graphon analogue $W$ we note that the 
copula $C$ represents a distribution on the
unit square and thus may be written as a composition of the distribution function of the first 
variable and the conditional distribution function $G(x,\cdot)$ of the second variable, 
given the value $x$ for the first. Let 
\begin{equation*}
             W(x,y):=\inf\{z\in [0,1]:\, G(x,z)\ge y\}, \quad 0\le y\le 1,
\end{equation*}
be the corresponding quantile function. Further, let $(U_i)_{i\in\bN}$ and $(V_i)_{i\in\bN}$ 
be two independent sequences of independent $\unif(0,1)$-distributed random variables. 
By construction, $(X_i,Y_i)\eqdistr (U_i,W(U_i,V_i))$ for all $i\in\bN$, and by independence,
this even holds for the two sequences. Thus the random binary matrix 
$\tilde Z=(\tilde Z_{ij})_{i,j\in\bN}$ given by 
\begin{equation*}
    \tilde Z_{ij}=\begin{cases}
                 1, &\text{if }i<j\text{ and } W(U_i,V_i)<W(U_j,V_j),\\
                 1, &\text{if }i>j\text{ and } U_i<U_j,\\
                 0, &\text{else,}  
            \end{cases}
\end{equation*}
has the same distribution as $Z$. Of course, the conditional construction could have been 
done in the other direction. It is indeed quite common in these constructions that uniqueness of
the representation only holds up to some equivalence, which may be difficult to describe.

It is of interest to know whether an exchangeable 
array is ergodic. For exchangeable sequences $(X_i)_{i\in\bN}$ this means that 
$M_\infty\equiv \mu$ for some fixed distribution $\mu$, so that the $X_i$'s are
independent. In the permuton situation such extreme points correspond to a fixed 
copula, and this is always the case for the models in Theorem~\ref{thm:Hoppen}\,(c). 
One of the interesting aspects 
of the results in \cite{Bass} is the fact that, for separable permutations, the limit object 
is not degenerate. 

Above we have used random matrices in order to point out the similarity to random
graphs and graphons. A more abstract and more general approach has been developed 
in~\cite{GerDiss}, where cartesian products of total orders on $\bN$ are considered
and ergodicity is discussed in some detail.

\subsection{Some connections to statistical concepts}\label{subsec:suff}
 Markov chain boundaries and exchangeable sequences both lead
to parametrized families of probability distributions, where the parameter is the value of
$X_\infty$ in the first case and the value of $M_\infty$ in the second.
A relationship between the Martin boundary of Markov chains and parametric families has
been pointed out early by~\cite{Abra}, together with a connection to exponential families
which, roughly, appear as $h$-transforms in certain situations. Section 18 in~\cite{AldousSF} deals with sufficiency and mixtures in the context of exchangeable random structures. 
In the framework of
combinatorial Markov chains $(X_n)_{n\in\bN}$ discussed above, where we regard the value 
of $X_\infty$ in the Doob-Martin compactification as a parameter $\theta$ and the values 
of the first $n$ variables $X_1,\ldots,X_n$ as data, it is clear that $X_n$ is a sufficient
statistic for $\theta$ as the conditional law of the data given $X_n$ can be reconstructed 
from the cotransitions, and these do not depend on $\theta$. 
\cite{Lau74} used this as a basis for his concept of \emph{total sufficiency}; the paper 
also discusses a variety of related aspects. \cite{Lau88} emphasizes the role of the Martin
boundary as a basic concept.

 Further, both the set of distributions on
the boundary and the set of directing measures are convex, and `pure parameters' or
`extreme models' correspond to the ergodic case, which are extremal elements of 
the respective convex set. This aspect has appeared repeatedly in the above examples. 
\cite{DynkinSuff} introduced \emph{H-sufficiency}, dealing 
with the question of whether the convex set is a simplex (in a barycentric sense). The latter
points to a connection with geometry; see~\cite{BaGrEJS} for a recent contribution. 

Two final comments may  be in order, both are somewhat loose.
First, as seen above, large 
random discrete structures can be investigated through the Doob-Martin boundary of Markov
chains, or through the representation of exchangeable distributions. The first is based on 
the topological approach of constructing a limit for a sequence of growing objects, 
the second is based on the construction of an `asymptotic template' from which the 
sequence can be obtained by sampling, and is of a more measure-theoretical nature. (Poisson boundaries, which we did not discuss, may be seen as a measure-theoretical variant 
of the former.)  A similar distinction appears in connection with convex sets, with 
Choquet's theorem for topological vector spaces, and the barycenter approach 
motivated by potential theory.
As a second comment, we note that permutations appear naturally in classical nonparametric
statistics, where ranks are a basic tool, but that other combinatorial families may similarly be 
analyzed with a view towards statistical applications. Often the results on the discrete
structures themselves can serve as a basis for obtaining asymptotics for a variety of specific
aspects, in some resemblance to using functional limit theorems to obtain distributional
limit theorems for a variety of specific statistics.   

\section*{Acknowledgments} 
Over the years, discussions with Ludwig Baringhaus, Steve Evans, 
Julian Gerstenberg, Klaas Hagemann, Anton Wakolbinger and Wolfgang Woess 
were instrumental for my understanding of many of the above concepts.     
I am also grateful to the referees for their helpful and encouraging comments.

\section*{Conflict of interest}

The author declares that he has no conflict of interest.

\bibliographystyle{spbasic} 

\begin{thebibliography}{}

\bibitem[Abrahamse(1970)]{Abra}
Abrahamse, A.F. (1970)
A comparison between the Martin boundary theory and the theory of likelihood ratios.
Ann. Math. Statistics 41: 1064--1067.

\bibitem[Aldous(1985)]{AldousSF} 
Aldous, D.J. (1985)
Exchangeability and related topics. 
Lecture Notes in Math. 1117:  1–-198. Springer, Berlin.

\bibitem[Arratia et al(2003)]{ABT}
Arratia, R., Barbour, A.D., Tavar\'{e}, S. (2003)
Logarithmic combinatorial structures: a probabilistic approach.
European Mathematical Society, Z\"{u}rich.

\bibitem[Austin(2008)]{Austin} 
Austin, T. (2008)
On exchangeable random variables and the statistics of large graphs and hypergraphs. 
Probability Surveys 5: 80--145.

\bibitem[Austin(2013)]{AustinLN}
Austin, T. (2013)
Exchangeable random arrays.
Lecture Notes, 
\texttt{www.math.ucla.edu/$\sim$tim/notes.html}.

\bibitem[Baringhaus and Gr\"ubel(2021)]{BaGrEJS}
Baringhaus, L., Gr\"ubel, R. (2021)
Discrete mixture representations of parametric distribution families: Geometry and statistics.
Electronic Journal of Statistics 15: 37-–70

\bibitem[Baringhaus and Gr\"ubel(2022)]{BaGrDeldat}
Baringhaus, L., Gr\"ubel, R. (2022)
Random permutations generated by delay models and estimation of delay distributions.
{\em In preparation.} 

\bibitem[Bassino et al(2018)]{Bass} 
Bassino, F., Bouvel, M., Féray, V., Gerin, L., Pierrot, A. (2018) 
The Brownian limit of separable permutations. 
Ann. Probab. 46: 2134--2189. 

\bibitem[Bona(2004)]{Bona} 
Bona, M. (2004) 
Combinatorics of Permutations.
Chapman and Hall, Boca Raton.

\bibitem[Ceccherini-Silberstein et al(2008)]{CSHarm}
Ceccherini-Silberstein, T., Scarabotti, F., Tolli, F. (2008)
Harmonic analysis on finite groups.
Cambridge University Press, Cambridge.

\bibitem[Cox and Smith(1961)]{CoxSmith}
Cox, D.R., Smith, W.L. (1961) 
Queues. Chapman and Hall, London.

\bibitem[Crudele et al(2023)]{sixperm}
Crudele, G., Dukes, P. and Noel, J.A. (2023)
Six permutation patterns force quasirandomness.
Available at: {\tt arXiv:2303.04776}.

\bibitem[Diaconis(1988)]{Diaconis} 
Diaconis, P. (1988). Group representations in probability and statistics. 
Institute of Mathematical Statistics Lecture Notes---Monograph Series 11. Hayward, CA.

\bibitem[Diaconis and Janson(1988)]{DiaconisJanson} 
Diaconis, P., Janson, S. (2008)
Graph limits and exchangeable random graphs.
Rend. Mat. Appl. 28: 33--61.

\bibitem[Doob(1959)]{Doob} 
Doob, J.L. (1959) 
Discrete potential theory and boundaries. 
J. Math. Mech. 8: 433--458.

\bibitem[Dynkin(1978)]{DynkinSuff}
Dynkin, E.B. (1978). Sufficient statistics and extreme points. 
Ann. Probab. 6: 705--730.

\bibitem[Evans al(2012)]{EGW1}
Evans, S.N., Gr{\"u}bel, R., Wakolbinger, A. (2012)
Trickle-down processes and their boundaries.
Electron. J. Probab. 17: 58pp.

\bibitem[Evans et al(2017)]{EGW2}
Evans, S.N., Gr{\"u}bel, R., Wakolbinger, A. (2017)
Doob--{M}artin boundary of {R}{\'e}my's tree growth chain.
Ann. Probab. 45: 225-277.

\bibitem[Feller(1968)]{FellerI}  
Feller, W. (1968) 
An Introduction to Probability Theory and Its Applications. Vol. I, 3rd ed..
Wiley, New York.

\bibitem[Gerstenberg(2018)]{GerDiss}
Gerstenberg, J. (2018)
Austauschbarkeit in Diskreten Strukturen: Simplizes und Filtrationen.
Dissertation. Leibniz Universit\"at Hannover.

\bibitem[Gerstenberg et al(2016)]{GGH}
Gerstenberg, J., Gr{\"u}bel, R., Hagemann, K. (2016)
A boundary theory approach to de Finetti's theorem.
Available at: https://arxiv.org/abs/1610.02561. 

\bibitem[Gr\"ubel(2013)]{GrSemBerKMK} 
Gr\"ubel, R. (2013)
Kombinatorische Markov-Ketten.
Math. Semesterber. 60: 185–-215

\bibitem[Gr\"ubel(2015)]{GrDMTCS} 
Gr\"ubel, R. (2015)
Persisting randomness in randomly growing discrete structures:
graphs and search trees.
Discrete Mathematics \& Theoretical Computer Science 18: 23pp.

\bibitem[Hagemann(2016)]{Hage}
Hagemann, K. (2016)
Doob-Martin-Theorie diskreter Markov-Ketten: Struktur und Anwendungen.
PhD thesis, Leibniz Universit\"at Hannover.

\bibitem[Hoppen et al(2013)]{Hopp} 
Hoppen, C., Kohayakawa, Y., Moreira, C.G., R\'ath, B., Sampaio, R.M. (2013)
Limits of permutation sequences.
J. Comb. Theory, Ser. B 103: 93--113.

\bibitem[Janson(2020)]{JansonAlg}
Janson, S. (2020) Patterns in random permutations avoiding some sets
of multiple patterns. Algorithmica 82, 616--641.

\bibitem[Kallenberg(2005)]{KallSym}
Kallenberg, O. (2005) Probabilistic symmetries and invariance principles.
Springer, New York.

\bibitem[Kelley(1955)]{Kelley} 
Kelley, J.L. (1955) 
General topology. Van Nostrand Company, Toronto-New York-London.

\bibitem[Kerov et al(1998)]{KOO}
Kerov, S.V., Okounkov, A., Olshanski, G. (1998)
The boundary of the Young graph with Jack edge multiplicities.
Int. Mathematics Research Notices: 173--199.

\bibitem[Kerov(2003)]{Kerov}
Kerov, S. V. (2003)
Asymptotic representation theory of the symmetric group and its applications in analysis.
Translations of Mathematical Monographs 219. American Mathematical Society, Providence, RI.

\bibitem[Kr\'al and Pikhurko(2013)]{Kral}
Kr{\'a}l, D., Pikhurko, O. (2013)
Quasirandom permutations are characterized by 4-point densities.
Geom. Funct. Anal. Vol. 23: 570–-579.

\bibitem[Lauritzen(1974)]{Lau74}
Lauritzen, S. L. (1974)
Sufficiency, prediction and extreme models.
Scand. J. Statist. 1: 128--134.

\bibitem[Lauritzen(1988)]{Lau88}
Lauritzen, S.L. (1988)
Extremal families and systems of sufficient statistics.
Lecture Notes in Statistics 49. Springer, Berlin.

\bibitem[Lovasz(2012)]{Lov} 
Lov\'asz, L. (2012) 
Large networks and graph limits. American Mathematical Society Colloquium
Publications 60. Amer. Math. Soc., Providence, RI.

\bibitem[M\'{e}liot(2017)]{Mel}
M\'{e}liot, P.-L. (2017)
Representation theory of symmetric groups.
CRC Press, Boca Raton.

\bibitem[Nelson(2006)]{Nelson} 
Nelson, R.B. (2006) 
An introduction to copulas, 2nd ed. Springer, New York.

\bibitem[Romik(2015)]{Romik}
Romik, D. (2015)
The surprising mathematics of longest increasing subsequences.
Institute of Mathematical Statistics Textbooks 4,
Cambridge University Press, New York.

\bibitem[Shi et al(2022)]{Drton}
Shi, H., Drton, M. and Han, F. (2022)
On the power of {Chatterjee}'s rank correlation.
Biometrika 109, 317--333.

\bibitem[van der Vaart(1998)]{vdV} 
van der Vaart, A.W. (1998)
Asymptotic Statistics. Cambridge University Press, Cambridge.

\bibitem[Woess(2009)]{WoessDMC}
Woess, W. (2009)
Denumerable Markov chains. 
European Mathematical Society, Z\"urich.

\bibitem[Yanagimoto(1970)]{Yana}
Yanagimoto, T. (1970) On measures of association and a related problem.
Ann. Inst. Stat. Math. 22, 57--63.

	
\end{thebibliography}

{}

\end{document}